# Entropy of convolutions on the circle

By Elon Lindenstrauss, David Meiri, and Yuval Peres

## Abstract

Given ergodic $p$-invariant measures $\{\mu_i\}$ on the 1-torus $\mathbb{T} = \mathbb{R}/\mathbb{Z}$, we give a sharp condition on their entropies, guaranteeing that the entropy of the convolution $\mu_1 * \cdots * \mu_n$ converges to $\log p$. We also prove a variant of this result for joinings of full entropy on $\mathbb{T}^\mathbb{N}$. In conjunction with a method of Host, this yields the following. Denote $\sigma_q(x) = qx \pmod 1$. Then for every $p$-invariant ergodic $\mu$ with positive entropy, $\frac{1}{N}\sum_{n=0}^{N-1} \sigma_{c_n}\mu$ converges weak* to Lebesgue measure as $N \longrightarrow \infty$, under a certain mild combinatorial condition on $\{c_k\}$. (For instance, the condition is satisfied if $p = 10$ and $c_k = 2^k + 6^k$ or $c_k = 2^{2^k}$.) This extends a result of Johnson and Rudolph, who considered the sequence $c_k = q^k$ when $p$ and $q$ are multiplicatively independent.

We also obtain the following corollary concerning Hausdorff dimension of sum sets: For any sequence $\{S_i\}$ of $p$-invariant closed subsets of $\mathbb{T}$, if $\sum \dim_H(S_i)/|\log \dim_H(S_i)| = \infty$, then $\dim_H(S_1 + \cdots + S_n) \longrightarrow 1$.

## 1. Introduction

Let $p \geq 2$ be any integer ($p$ need not be prime), and $\mathbb{T} = \mathbb{R}/\mathbb{Z}$ the 1-torus. Denote by $\sigma_p$ the $p$-to-one map $x \mapsto px \pmod 1$. The pair $(\mathbb{T}, \sigma_p)$ is a dynamical system that has additional structure: $\mathbb{T}$ is a commutative group (with the group operation being addition mod 1), and $\sigma_p$ is an endomorphism of it. Even in such a simple system, the interaction between the dynamics and the algebraic structure of $\mathbb{T}$ can be quite subtle; the present work continues the study of this interaction, inspired by the fundamental work of Furstenberg [8].

Say that a measure $\mu$ on $\mathbb{T}$ is $p$-invariant if $\sigma_p\mu = \mu$, where for every set $A \subset \mathbb{T}$
$$(\sigma_p\mu)(A) \stackrel{\text{def}}{=} \mu(\sigma_p^{-1}A).$$

(All measures we consider are Borel probability measures.) Lebesgue measure on $\mathbb{T}$, denoted $\lambda$, has entropy $\log p$ with respect to $\sigma_p$, and is the unique $p$-invariant measure of maximal entropy. Given two $p$-invariant measures $\mu$ and $\nu$, the group structure of $\mathbb{T}$ naturally yields another $p$-invariant measure — the convolution $\mu * \nu$.

Our main results, Theorems 1.1 and 1.8, concern the entropy growth for convolutions of $p$-invariant measures and their ergodic components. These results have applications to the Hausdorff dimension of sum sets and to genericity of the orbits of measures with positive entropy under multiplication by certain integer sequences.

THEOREM 1.1 (the Convolution Theorem). *Let $\{\mu_i\}$ be a sequence of $p$-invariant and ergodic measures on $\mathbb{T}$ whose normalized base-$p$ entropies $h_i = h(\mu_i, \sigma_p)/\log p$ satisfy*

$$(1) \qquad \sum_{i=1}^{\infty} \frac{h_i}{|\log h_i|} = \infty.$$

*Then*

$$h(\mu_1 * \cdots * \mu_n, \sigma_p) \longrightarrow \log p \quad \text{monotonically,} \quad \text{as } n \longrightarrow \infty.$$

*In particular, $\mu_1 * \cdots * \mu_n \longrightarrow \lambda$ weak$^*$ and in the $\bar{d}$ metric (with respect to the base-$p$ partition).*

It is relatively easy to see that under hypotheses of the theorem, $\mu_1 * \cdots * \mu_n \longrightarrow \lambda$ weak$^*$. This means that

$$\int f(x)\, d\mu_1 * \cdots * \mu_n \longrightarrow \int f(x)\, d\lambda$$

for all continuous $f$, and gives very little information on the dynamics of $\mu_1 * \cdots * \mu_n$.

The convergence of the entropy to $\log p$ is equivalent to the much stronger statement that $\mu_1 * \cdots * \mu_n \longrightarrow \lambda$ in the $\bar{d}$ metric. As we will not use this metric in our arguments, we only recall its definition and refer the reader to Rudolph [17] for further information. Consider two $p$-invariant measures $\nu_1$ and $\nu_2$ on $\mathbb{T}$. What $\bar{d}(\nu_1, \nu_2) < \epsilon$ means is that there exists a $p$-invariant measure $\tilde{\nu}$ on $\mathbb{T}^2$ that projects to $\nu_1$ in the first coordinate and to $\nu_2$ in the second coordinate, such that for $\tilde{\nu}$-almost every $(x, y) \in \mathbb{T}^2$, the set of integers $k \geq 1$ for which the $k^{\text{th}}$ digits in the base $p$ expansions of $x$ and $y$ differ, has asymptotic density less than $\epsilon$. Once we establish the entropy convergence, the $\bar{d}$ convergence of $\mu_1 * \cdots * \mu_n$ to $\lambda$ is an immediate corollary of the fact that $\lambda$ is a Bernoulli measure, and hence is finitely determined; see Rudolph [17, Chaps. 6 and 7], for the relevant definitions and proofs.

The entropy condition in Theorem 1.1 is sharp: if $\{h_i\}$ is a sequence of numbers in the range $(0, 1)$ with $\sum h_i/|\log h_i| < \infty$, then there exists a sequence of $p$-invariant ergodic measures $\{\mu_i\}$, such that $h(\mu_i, \sigma_p)/\log p = h_i$, yet $\mu_1 * \cdots * \mu_n$ does not converge to Lebesgue measure $\lambda$ even in the weak* topology; see Example 10.2.

The convolution theorem has implications for Hausdorff dimension of sum sets:

COROLLARY 1.2. *Let $\{S_i\}$ be a sequence of $p$-invariant closed subsets of $\mathbb{T}$, and suppose that*
$$\sum_{i=1}^{\infty} \frac{\dim_H(S_i)}{|\log \dim_H(S_i)|} = \infty.$$
*Then $\dim_H(S_1 + \cdots + S_n) \longrightarrow 1$.*

By Furstenberg [8, III, 2], the conclusion is equivalent to $h_{\text{top}}(S_1 + \cdots + S_n, \sigma_p) \to \log p$.

If the measures $\mu_i$ are not weakly-mixing, the measure $\mu_1 * \cdots * \mu_n$ might be nonergodic, with different fibers carrying different entropies (see Example 9.4). The reason for this is that for two measures $\mu$ and $\nu$ on $\mathbb{T}$, the convolution $\mu * \nu$ is the projection of the product measure $\mu \times \nu$ on $\mathbb{T} \times \mathbb{T}$ to $\mathbb{T}$; if $\mu$ and $\nu$ are not weakly mixing, then $\mu \times \nu$ need not be ergodic (indeed, $\mu$ is weakly mixing if and only if $\mu \times \mu$ is ergodic). In this case, ergodic components of $\mu \times \nu$ can project to ergodic components of $\mu * \nu$ with distinct entropies (however, it is easy to see that the entropy of almost all ergodic components of $\mu * \nu$ is at least the entropy of $\mu$ — see Corollary 9.3). In this more general situation, when $\mu_i$ are not assumed to be weakly mixing, it turns out to be both natural and important for the applications to give more accurate information than that provided by Theorem 1.1 regarding the ergodic components of $\mu_1 * \cdots * \mu_n$. This is done in Theorem 1.8 below; the proof of this more detailed result is rather delicate. As before, such a result can be used to obtain an estimate of Hausdorff dimension of sum sets:

THEOREM 1.3. *Let $\{\mu_i\}$ be a sequence of $p$-invariant ergodic measures on $\mathbb{T}$, and suppose that $\inf_i h(\mu_i, \sigma_p) > 0$. Let $\{S_i\}$ be a sequence of Borel subsets of $\mathbb{T}$, and suppose that $\mu_i(S_i) > 0$ for all $i \geq 1$. Then $\dim_H(S_1 + \cdots + S_n) \longrightarrow 1$.*

Note that the the measures $\mu_i(S_i)$ can tend to zero arbitrarily fast.

Our initial motivation for studying entropy of convolutions of $p$-invariant measures was to find conditions on a sequence of integers $\{c_n\}$ and a measure $\mu$, which imply that $\mu$ is $\{c_n\}$-*generic*, i.e., the averages $\frac{1}{N}\sum_{n=0}^{N-1} \sigma_{c_n}\mu$ converge weakly to Lebesgue measure $\lambda$. In certain cases, this could be established for ergodic $p$-invariant measures of sufficiently high entropy, and the idea was to deduce $\{c_n\}$-genericity for all ergodic $p$-invariant measures of positive entropy by repeated convolutions. Indeed, we establish Theorem 1.4 below using this scheme.

The works of Host [9] and Meiri [13] indicate that combinatorial properties of $\{c_n\}$ can be used to prove $\{c_n\}$-genericity for invariant measures with positive entropy. The combinatorial property we need is weaker than those assumed in the quoted papers; we require that the number of pairs among the first $p^n$ elements of $\{c_n\}$ which are congruent mod $p^n$, is exponentially smaller than $p^{2n}$.

*Definition* 1.1. Given an integer-valued sequence $\{c_n\}$ and an integer $p > 1$, we define the *p-adic collision exponent* of the sequence as

$$\Gamma_p(\{c_n\}) = \limsup_{n \to \infty} \frac{\log |\{0 \leq k, \ell < p^n : c_k \equiv c_\ell \pmod{p^n}\}|}{n \log p}.$$

Since pairs with $k = \ell$ are allowed, we always have $1 \leq \Gamma_p(\{c_n\}) \leq 2$. For example, whenever $p, q > 1$ are relatively prime, the $p$-adic collision exponent of $\{q^n\}$ is 1. If we assume only that there is some prime factor of $p$ which does not divide $q$, then $\{q^n\}$ has $p$-adic collision exponent $< 2$. (See §3.2 for more details and refinements of this definition.)

Using this definition we can state our results on genericity of orbits of $p$-invariant measures:

THEOREM 1.4. *Let $\{c_n\}$ be a sequence with $p$-adic collision exponent $< 2$, for some $p > 1$. Then any $p$-invariant ergodic measure $\mu$ on $\mathbb{T}$ with positive entropy is $\{c_n\}$-generic, i.e.,*

$$\text{(2)} \qquad \frac{1}{N} \sum_{n=0}^{N-1} \sigma_{c_n} \mu \longrightarrow \lambda \quad \text{in the weak* topology.}$$

For instance, the hypothesis is satisfied if $p = 10$ and $c_n = 2^n + 6^n$ or $c_n = 2^{2^n}$; see Section 4 for other examples. In fact, the following stronger form of convergence holds. Recall that the space of probability measures on $\mathbb{T}$ endowed with the weak* topology is a compact metric space, and take $\rho_*$ be some metric on it.

THEOREM 1.5. *Under the conditions of Theorem 1.4, $\mu$ is $\{c_n\}$-normal in probability, i.e.,*

$$\text{(3)} \qquad \int \rho_*\Big(\frac{1}{N} \sum_{n=0}^{N-1} \delta_{c_n x}, \lambda\Big) d\mu(x) \longrightarrow 0.$$

To be more specific, let

$$\hat{\mu}(k) = \int e^{2\pi i k x} d\mu .$$

Then $\rho_*(\mu, \nu) = \sum_{k=-\infty}^{\infty} 2^{-|k|} |\hat{\mu}(k) - \hat{\nu}(k)|^2$ is a metric on the space of all $p$-invariant probability measures on $\mathbb{T}$ with the weak* topology. Define for any integer $k \neq 0$,

$$g_N^{(k)}(x) \stackrel{\text{def}}{=} \frac{1}{N} \sum_{n=0}^{N-1} e(kc_n x), \tag{4}$$

where $e(x) \stackrel{\text{def}}{=} e^{2\pi i x}$. Then (2) is equivalent to for all $k \neq 0$, $\int g_N^{(k)}(x)\, d\mu \longrightarrow 0$, while (3) is equivalent to the stronger property: for all $k \neq 0$, $\int |g_N^{(k)}(x)|^2\, d\mu \longrightarrow 0$.

The case $c_n = q^n$ is known, and inspired our general investigation. Even though there are multiplicatively-independent $p$ and $q$ with $\Gamma_p(\{q^n\}) = 2$, the following is still a corollary of the above results.

COROLLARY 1.6 (Johnson and Rudolph [11, Thm. 8.6]). *Suppose that $p, q > 1$ are multiplicatively-independent and $\mu$ is a $p$-invariant and ergodic measure on $\mathbb{T}$ with positive entropy. Then $\mu$ is $\{q^n\}$-normal in probability.*

Our proofs of Theorems 1.4–1.5 use two main tools. Host [9] developed a harmonic analysis method which is most powerful when the entropy of $\mu$ is large. The following general result then allows us to use Host's method for all measures with nonzero entropy by reduction (via convolutions) to the case where the entropy of the measure $\mu$ is sufficiently high.

THEOREM 1.7 (the Bootstrap Lemma). *Suppose that $\mathcal{C}$ is a class of $p$-invariant measures on $\mathbb{T}$ with the following properties:*

(i) *If $\mu$ is $p$-invariant and ergodic and $\mu * \mu \in \mathcal{C}$, then $\mu \in \mathcal{C}$.*

(ii) *If $\mu$ is $p$-invariant and almost every ergodic component of $\mu$ is in $\mathcal{C}$, then $\mu \in \mathcal{C}$.*

(iii) *There exists some constant $h_0 < \log p$ such that every $p$-invariant and ergodic measure $\mu$ with $h(\mu, \sigma_p) > h_0$ is in $\mathcal{C}$.*

*Then $\mathcal{C}$ contains all $p$-invariant ergodic measures with positive entropy.*

We derive Theorem 1.7 from a variant of Theorem 1.1 for joinings of full entropy:

*Definition* 1.2. Let $\{\mu_i\}_{i \geq 1}$ be a sequence of $p$-invariant and ergodic measures on $\mathbb{T}$. A measure $\nu^{(n)}$ on $\mathbb{T}^n$ is called a *joining* of $\mu_1, \ldots, \mu_n$ if

(i) $\nu^{(n)}$ is $\sigma_p \times \cdots \times \sigma_p$ invariant,

(ii) The projection of $\nu^{(n)}$ on the $i^{\text{th}}$ coordinate is $\mu_i$, for $i = 1, \ldots, n$.

The measure $\nu^{(n)}$ is called a *joining of full entropy* of $\mu_1, \ldots, \mu_n$ if in addition

$$h(\nu^{(n)}, \sigma_p \times \cdots \times \sigma_p) = \sum_{i=1}^{n} h(\mu_i, \sigma_p).$$

A measure $\tilde{\mu}$ on $\mathbb{T}^{\mathbb{N}}$ is called a joining of full entropy of $\{\mu_i\}_{i=1}^{\infty}$ if for every $n$, the projection of $\tilde{\mu}$ to the first $n$ coordinates is a joining of full entropy of $\mu_1, \ldots, \mu_n$.

THEOREM 1.8. *Let $\{\mu_i\}_{i=1}^{\infty}$ be a sequence of $p$-invariant and ergodic measures on $\mathbb{T}$ such that $\inf_i h(\mu_i, \sigma_p) > 0$. Suppose that $\tilde{\mu}$ is a joining of full entropy of $\{\mu_i\}$. Define $\Theta^n : \mathbb{T}^{\mathbb{N}} \to \mathbb{T}$ by $\Theta^n(x) = x_1 + \cdots + x_n \pmod{1}$. Then*

$$h(\Theta^n \tilde{\mu}, \sigma_p) \longrightarrow \log p \quad \text{monotonically}, \quad \text{as } n \longrightarrow \infty.$$

Theorem 1.8 is not valid under the weaker entropy assumptions of Theorem 1.1. Indeed, it is possible to find a joining of full entropy $\tilde{\mu}$ with entropies satisfying (1), such that $\Theta^n \tilde{\mu}$ does not even converge weak* to $\lambda$. See Example 10.5.

1.1. *Background.* In Furstenberg [8] many aspects of the dynamics of $(\mathbb{T}, \sigma_p)$ are discussed, and in particular it is shown that there are no nontrivial $\sigma_p$-invariant closed subsets of $\mathbb{T}$ that are also invariant under $\sigma_q$ for $p$ and $q$ multiplicatively independent (i.e. $\log p / \log q \notin \mathbb{Q}$), the trivial examples being the whole of $\mathbb{T}$ and some finite sets of periodic points. Furstenberg conjectured the following stronger result:

FURSTENBERG'S CONJECTURE. *The only ergodic invariant measures for the semi-group of circle endomorphisms generated by $\sigma_p$, and $\sigma_q$ for $p$ and $q$ multiplicatively independent are Lebesgue measure $\lambda$, and atomic measures concentrated on periodic orbits.*

Most of the research on the dynamics of $(\mathbb{T}, \sigma_p)$ has been related to this conjecture. It has been proven for measures with positive entropy — a partial result was proved by Lyons [12], and the case of $p$ and $q$ relatively prime was settled by Rudolph [18]. The case of $p$ and $q$ multiplicatively independent but not relatively prime is harder, and was proved by Johnson [10]. Another argument along the lines of Lyons [12] for the multiplicatively independent case was given by Feldman [6]. To tackle the case of measures with zero entropy it seems that a totally different method is needed.

In Feldman and Smorodinsky [7], Johnson and Rudolph [11], and Host [9], the measure is only assumed to be $\sigma_p$-invariant, and the authors consider the action of $\sigma_{c_n}$ on $\mu$ for the special case of $c_n = q^n$. As shown in Meiri [13], the methods of Host imply the following result:

THEOREM A.    *For $p > 1$, let $\{c_n\}$ be a sequence with p-adic collision exponent 1. Then any p-invariant ergodic measure $\mu$ with positive entropy is $\{c_n\}$-normal a.e., i.e., $\{c_n x \pmod 1\}$ is uniformly distributed for $\mu$-almost every $x \in \mathbb{T}$.*

This theorem gives a significantly stronger statement than Theorems 1.4–1.5, but for a smaller class of sequences $\{c_n\}$ — and in particular gives no information for the case $c_n = q^n$, with $p$ and $q$ not relatively prime. (Recalling the definition of $g_N^{(k)}$ in (4), this assertion is equivalent to for all $k \neq 0$, $g_N^{(k)} \longrightarrow 0$ $\mu$–a.e.) A more detailed history of this problem can be found in Host [9] and Meiri [15].

For the convolution results (Thms. 1.1 and 1.8), the case of measures with zero entropy is not interesting, as the entropy of a convolution of measures, each having zero entropy, has zero entropy, and hence does not converge $\bar{d}$ to $\lambda$.

Regarding Theorems 1.2–1.3, one can ask if a stronger conclusion holds, namely, that $S_1 + \cdots + S_n$ contains an interval for all $n$ sufficiently large. Brown and Williamson [4] showed that if $\mu$ is a measure on $\mathbb{T}$ which makes the digits in base $p$ i.i.d. nondegenerate variables, and $\mu(S_i) > 0$ for all $i$, then this stronger assertion is true. However, under our weaker assumptions, this conclusion is not valid, even if all the sets $S_i$ coincide: By Furstenberg [8, Thm. III.2] there exists a minimal $p$-invariant closed set $S \subset \mathbb{T}$ with positive dimension. By Proposition IV.1 of the same article, the sum sets $S + \cdots + S$ of any finite order have Lebesgue measure zero. For more information about Lebesgue measure of sum sets see Brown, Keane, Moran, and Pearce [3].

1.2. *Overview.* The paper is organized as follows. In Section 2 we show how one can deduce the Bootstrap Lemma from Theorem 1.8, and proceed to prove Theorem 1.4. In Section 3 we discuss convergence in probability and prove Theorem 1.5. In Section 4 we derive Corollary 1.6 from Theorem 1.4, and discuss $p$-adic combinatorial properties, and in particular we give an algorithm for computing the $p$-adic collision exponent of linear recursive sequences.

In Sections 5–7 we prove our main results, Theorems 1.1 and 1.8. The simplest case is proving Theorem 1.1 for prime $p$. In this case, however, one does not need to use a substantial part of the ideas behind the proof. We recommend for first reading to have in mind the case $p = 10$, with $\mu_i = \mu$ for all $i$ for Theorem 1.1. Theorem 1.8 is interesting already when $p$ is a prime (and we again recommend considering first the case where all $\mu_i$ are identical).

Section 5 contains results about finite cyclic groups which are crucial to the proof of the convolution theorems. Lemmas 5.1 and 5.2 study convolutions of measures on a finite cyclic group and contain one key idea in the proof, namely that the convolution of a sequence of measures on a finite cyclic group of order $N$ (we shall use $N = p^k$) tends to be invariant under a subgroup that will typically be rather large (in the cases we will be interested in, this subgroup will be of order approximately $p^{\alpha k}$ for some $0 < \alpha < 1$). Lemma 5.3 shows that if a measure on $\mathbb{Z}/p^k\mathbb{Z}$ is almost invariant under a subgroup of size $p^{\alpha k}$, the distribution of the $\alpha k$ high order digits is nearly uniform.

In Section 6 we begin to show how convolutions of measures on $\mathbb{Z}/p^k\mathbb{Z}$ relate to convolutions of measures on $\mathbb{T}$, where we get measures on $\mathbb{Z}/p^k\mathbb{Z}$ from measures on $\mathbb{T}$ by considering the conditional distribution of the first $k$ digits in the base-$p$ expansion of $x \in \mathbb{T}$, given the rest of the digits. In Section 7 we continue in this approach and prove Theorem 1.1. The basic argument is that if the entropy of $\mu_1 * \cdots * \mu_n$ is almost

$$\sup_{N \in \mathbb{N}} h(\mu_1 * \ldots * \mu_N, \sigma_p),$$

then for any $k \geq 1$ the distribution of the first $k$ digits of $x$ given the rest of the digits ($x$ chosen according to $\mu_1 * \cdots * \mu_n$) must be nearly invariant under a subgroup $G \subset \mathbb{Z}/p^k\mathbb{Z}$ of size $p^{\alpha k}$ — for if it is not, then the entropy of $\mu_1 * \cdots * \mu_n$ can be significantly increased by further convolutions. This implies that the first $\alpha k$ digits of $x$ are distributed nearly uniformly. Since $k$ is arbitrary, it follows that $h(\mu_1 * \cdots * \mu_n, \sigma_p) \simeq \log p$.

In Section 8 we prove Theorem 1.8. The main observation is that if one chooses an element in $\mathbb{T}^\mathbb{N}$ according to a joining of full entropy, the components of this vector are almost independent (Lemma 8.1). This allows us to prove Theorem 1.8 along similar lines as Theorem 1.1. In Section 9 we use the connection between entropy of measures and Hausdorff dimension, to derive Theorems 1.2–1.3 from the results about convolutions of measures. Section 10 contains some concluding remarks, questions and examples.

## 2. Proof of the Bootstrap Lemma and Theorem 1.4

In order to derive the Bootstrap Lemma from Theorem 1.8, we need the following lemma:

LEMMA 2.1. *Let $\mu$ denote a p-invariant and ergodic measure on $\mathbb{T}$, and let $\mu \times \mu = \int \nu_\gamma \, d\gamma$ denote the ergodic decomposition of $\mu \times \mu$. Then $\nu_\gamma$ is a self-joining of $\mu$ for a.e. $\gamma$, and $h(\nu_\gamma, \sigma_p \times \sigma_p) = 2h(\mu, \sigma_p)$ a.e.*

*Proof.* By ergodicity of $\mu$, almost-surely $\nu_\gamma$ projects to $\mu$ in both coordinates. Obviously, $h(\nu_\gamma, \sigma_p \times \sigma_p) \leq h(\mu \times \mu, \sigma_p \times \sigma_p) = 2h(\mu, \sigma_p)$. On the other hand, by Rokhlin's theorem, $h(\mu \times \mu, \sigma_p \times \sigma_p) = \int h(\nu_\gamma, \sigma_p \times \sigma_p) \, d\gamma$, therefore $h(\nu_\gamma, \sigma_p \times \sigma_p) = 2h(\mu, \sigma_p)$ a.e. □

*Proof of Theorem* 1.7 (*the Bootstrap Lemma*). Let $\mu$ be some $p$-invariant and ergodic measure with positive entropy. We claim that $\mu \in \mathcal{C}$. Suppose that this is not the case, and write $\nu^{(1)} \stackrel{\text{def}}{=} \mu$ and $h \stackrel{\text{def}}{=} h(\mu, \sigma_p) > 0$. From property (i) of $\mathcal{C}$, also $\mu * \mu \notin \mathcal{C}$. Let $\mu \times \mu = \int \nu_\gamma \, d\gamma$ denote the ergodic decomposition of $\mu \times \mu$ with respect to $\sigma_p \times \sigma_p$. Then

$$\mu * \mu = \Theta^2(\mu \times \mu) = \int \Theta^2(\nu_\gamma) \, d\gamma,$$

where as before $\Theta^k(x_1, \ldots, x_k) \stackrel{\text{def}}{=} x_1 + \cdots + x_k \pmod 1$. From property (ii), for a set of $\gamma$ with positive measure, $\Theta^2(\nu_\gamma) \notin \mathcal{C}$. By Lemma 2.1, there exists an ergodic component of $\mu \times \mu$, which we designate $\nu^{(2)}$, such that $\nu^{(2)}$ is an ergodic self-joining of $\mu$ with entropy $2h$ and $\Theta^2(\nu^{(2)}) \notin \mathcal{C}$. Apply now the same procedure to $\nu^{(2)}$, finding an ergodic component $\nu^{(4)}$ of $\nu^{(2)} \times \nu^{(2)}$, such that $\nu^{(4)}$ is a self-joining of $\nu^{(2)}$ with entropy $4h$, and $\Theta^4(\nu^{(4)}) \notin \mathcal{C}$. Continuing this way we obtain a sequence of measures $\{\nu^{(k)}\}$, defined for $k$ a power of 2, such that $\Theta^k(\nu^{(k)}) \notin \mathcal{C}$, $\nu^{(k)}$ is an ergodic self-joining of $\mu$ in $\mathbb{T}^k$, and $h(\nu^{(k)}, \sigma_p \times \cdots \times \sigma_p) = kh$. Define $\nu^{(k)}$ for other values of $k$ by projection, and let $\tilde{\mu}$ be the inverse limit of these measures, defined on $\mathbb{T}^\mathbb{N}$. Then $\tilde{\mu}$ is a joining of full entropy. Applying Theorem 1.8 we conclude that $h(\Theta^k \nu^{(k)}, \sigma_p) = h(\Theta^k \tilde{\mu}, \sigma_p) \longrightarrow \log p$. As $\Theta^k \tilde{\mu} \notin \mathcal{C}$ for $k$ a power of 2, this contradicts property (iii). $\square$

We next deduce Theorem 1.4 from the Bootstrap Lemma. Given an integer-valued sequence $\{c_n\}$, say that a measure $\mu$ is $\{c_n\}$-*generic* if

$$\frac{1}{N} \sum_{n=0}^{N-1} \sigma_{c_n} \mu \longrightarrow \lambda \text{ weakly}, \quad \text{as } N \longrightarrow \infty.$$

The basic tool we use is the following observation in Meiri [13, Thm. 3.1, and note in §8, problem 3], based on Host [9]:

PROPOSITION 2.2. *Fix an integer $p > 1$ and a sequence $\{c_k\}$ with $p$-adic collision exponent $< 2$. Then there exists a constant $h_0 < \log p$ such that every $p$-invariant and ergodic measure $\mu$ with $h(\mu, \sigma_p) > h_0$ is $\{c_n\}$-normal a.e.; in particular, $\mu$ is $\{c_n\}$-generic.*

*Remark.* This proposition, as well as Theorems 1.4–1.5, hold under the weaker assumption that the *reduced* $p$-adic exponent of $\{c_k\}$ is less than 2; see Section 3.2 for the definition and more details.

The following lemma is a consequence of the extremality of ergodic measures:

LEMMA 2.3 (Johnson and Rudolph [11]). *Suppose that $\nu, \nu_1, \nu_2, \ldots$ are invariant measures, and that $\nu$ is ergodic. Suppose also that*

$$\frac{1}{N} \sum_{n=1}^{N} \nu_n \longrightarrow \nu \text{ weakly.}$$

*Then there exists a zero-density set $J \subset \mathbb{N}$ such that $\nu_n \longrightarrow \nu$ (weakly) as $n \underset{n \notin J}{\longrightarrow} \infty$.*

PROPOSITION 2.4. *An invariant measure $\mu$ is $\{c_n\}$-generic if and only if $\mu * \mu$ is $\{c_n\}$-generic.*

*Proof.* By the previous lemma, $\mu$ is $\{c_n\}$-generic if and only if there exists a zero-density set $J \subset \mathbb{N}$ such that $\sigma_{c_n} \mu \longrightarrow \lambda$ as $n \underset{n \notin J}{\longrightarrow} \infty$. This is equivalent to

(5) $$\lim_{n \notin J} \hat{\mu}(ac_n) = 0, \text{ for all } a \in \mathbb{Z} \smallsetminus \{0\}.$$

Since $\widehat{\mu * \mu} = \hat{\mu}^2$, equation (5) holds if and only if it holds when $\mu$ is replaced by $\mu * \mu$. $\square$

*Proof of Theorem 1.4.* From Propositions 2.2 and 2.4 it follows that the class of $p$-invariant measures which are $\{c_n\}$-generic satisfies the conditions of the Bootstrap Lemma (Theorem 1.7), and the assertion follows. $\square$

## 3. Convergence in probability and proof of Theorem 1.5

Recall the definition (4) of $g_N^{(k)}$. Given a measure $\mu$, define $\mu_\#$ by $\mu_\#(A) = \mu(-A)$.

LEMMA 3.1. *For any measure $\mu$,*

$$\int |g_N^{(k)}(x)|^2 \, d\mu \leq \left( \int |g_N^{(k)}(x)|^2 \, d\mu * \mu_\# \right)^{\frac{1}{2}}.$$

*Proof.* As $\widehat{\mu * \mu_\#} = |\hat{\mu}|^2$, by Cauchy-Schwarz we have

$$\int |g_N^{(k)}(x)|^2 \, d\mu = \frac{1}{N^2} \sum_{m,l=0}^{N-1} \int e(k(c_l - c_m)x) \, d\mu = \frac{1}{N^2} \sum_{m,l=0}^{N-1} \hat{\mu}(k(c_l - c_m))$$

$$\leq \frac{1}{N} \left( \sum_{m,l=0}^{N-1} |\hat{\mu}(k(c_l - c_m))|^2 \right)^{\frac{1}{2}} = \left( \int |g_N^{(k)}(x)|^2 \, d\mu * \mu_\# \right)^{\frac{1}{2}}.$$

$\square$

*Proof of Theorem* 1.5. A slight change in the proof of the Bootstrap Lemma lets us replace the first condition by the following:

(i') If $\mu$ is $p$-invariant and $\mu*\mu_\# \in \mathcal{C}$, then $\mu \in \mathcal{C}$. Taking $\mathcal{C}$ to be the class of $p$-invariant measures for which (3) holds, condition (i') is satisfied by the last lemma, (ii) is immediate by Lebesgue dominated convergence, while (iii) follows from Proposition 2.2. □

The type of convergence in Theorem 1.5 is, in general, stronger than weak convergence:

*Example.* Let $c_n = 2^{2^n}$. For $j = 0, 1$ consider random variables $X^{(j)} = \sum_{i=1}^{\infty} x_i^{(j)} 2^{-i}$, where $x_{k^2}^{(j)} = j$ for all $k$, and all other digits are i.i.d. uniform on $\{0,1\}$. Let $\mu_j$ be the distribution of $X_j$, and take $\mu = \frac{1}{2}(\mu_1 + \mu_2)$. Then $\mu$ is $\{2^{2^n}\}$-generic, but is not $\{2^{2^n}\}$-normal in probability.

However, there are cases where the stronger type of convergence follows from the weaker one. The following proposition was obtained by Johnson and Rudolph [11, §8] using general convex analysis. Here we give a more direct argument.

PROPOSITION 3.2. *Let $q > 1$. Then if $\mu$ is $\{q^n\}$-generic, it is also $\{q^n\}$-normal in probability.*

*Proof.* For any $L$, by Cauchy-Schwarz,

$$\left| \frac{1}{N} \sum_{n=0}^{N-1} g_L^{(k)} \circ \sigma_{q^n} \right|^2 \leq \frac{1}{N} \sum_{n=0}^{N-1} |g_L^{(k)}|^2 \circ \sigma_{q^n},$$

and so by (2) applied for $c_n = q^n$,

$$\limsup_{N\to\infty} \int \left| \frac{1}{N} \sum_{n=0}^{N-1} g_L^{(k)} \circ \sigma_{q^n} \right|^2 d\mu \leq \limsup_{N\to\infty} \frac{1}{N} \sum_{n=0}^{N-1} \int |g_L^{(k)}|^2 \circ \sigma_{q^n} d\mu$$

$$= \int |g_L^{(k)}|^2 d\lambda = \frac{1}{L}.$$

Since

$$\left| g_N^{(k)} - \frac{1}{N} \sum_{n=0}^{N-1} g_L^{(k)} \circ \sigma_{q^n} \right| \leq \frac{2L}{N},$$

we conclude that for any $L$,

$$\limsup_{N\to\infty} \int |g_N^{(k)}|^2 d\mu \leq \frac{1}{L}$$

it follows that $\mu$ is $\{q^n\}$-normal in probability. □

### 4. The $p$-adic collision exponent

Denote by $\mathrm{ord}_G(x)$ the order of an element $x$ in a finite group $G$. The following statement was noted by Host [9]:

PROPOSITION 4.1. *If $p, q > 1$ and $p, q$ are relatively prime, then there exists $\alpha > 0$ such that*

$$\text{for all } n \geq 1, \ \mathrm{ord}_{\mathbb{Z}/p^n\mathbb{Z}}(q) \geq \alpha p^n.$$

*In particular, in this case the $p$-adic collision exponent of $\{q^n\}$ is 1.*

PROPOSITION 4.2. *Suppose that $p, q > 1$ and there exists some prime factor $p_*$ of $p$ that does not divide $q$. Then the $p$-adic collision exponent of $\{q^n\}$ is $< 2$.*

*Proof.* Define $o_n \stackrel{\text{def}}{=} \mathrm{ord}_{\mathbb{Z}/p_*^n\mathbb{Z}}(q)$, and use Proposition 4.1 to find $\alpha > 0$ such that $o_n \geq \alpha p_*^n$. If $q^k \equiv q^\ell \pmod{p^n}$ then $p_*^n | q^{|k-\ell|} - 1$; hence $o_n | k - \ell$. Thus

$$\left|\{0 \leq k, \ell < p^n : q^k \equiv q^\ell \pmod{p^n}\}\right| \leq p^n \cdot (p^n / o_n).$$

Hence

$$\Gamma_p(\{q^n\}) = \limsup_{n \to \infty} \frac{\log\left|\{0 \leq k, \ell < p^n : q^k \equiv q^\ell \pmod{p^n}\}\right|}{n \log p}$$

$$\leq \limsup_{n \to \infty} \frac{\log(p^n \cdot (p^n / o_n))}{n \log p} \leq 2 - \frac{\log p_*}{\log p} < 2. \qquad \square$$

*Proof of Corollary* 1.6. By Proposition 3.2, it is enough to prove weak convergence. If some prime factor of $p$ does not divide $q$, we are done by the last proposition and Theorem 1.4. Otherwise, note the following simple observations:

(A) If the theorem holds replacing $q$ by some power $q^\ell$, then it holds for $q$ as well. The reason is that we can decompose

$$\frac{1}{N\ell} \sum_{n=0}^{N\ell-1} \sigma_{q^n} \mu = \frac{1}{\ell} \sum_{k=0}^{\ell-1} \frac{1}{N} \sum_{n=0}^{N-1} \sigma_{q^{\ell n}}(\sigma_{q^k} \mu)$$

and apply the theorem to the measures $\sigma_{q^k} \mu$, for $k = 0, \ldots, \ell - 1$.

(B) If $p | q$ and the theorem holds for $p$ and $q/p$, then it also holds for $p$ and $q$. This is because $\sigma_{q/p} \mu = \sigma_q \mu$ for $p$-invariant $\mu$.

For any multiplicatively independent $p, q$ we can find some $k, \ell$ such that $p^k | q^\ell$ and some prime factor of $p$ does not divide $q' = q^\ell / p^k$, and use the above. □

It is easy to see that Theorems 1.4–1.5 are not affected if we only assume that $\{c_n\}$ differs on a set of arbitrarily small density from a sequence with $p$-adic collision exponent smaller than 2. This motivates the following definition:

*Definition 4.1.* The *reduced $p$-adic collision exponent* of a sequence $\{c_n\}$ is

$$\Gamma'_p(\{c_n\}) = \lim_{\varepsilon \to 0} \inf_{\{c'_n\}} \Gamma_p(\{c'_n\}),$$

where $\{c'_n\}$ ranges over sequences which agree with $\{c_n\}$ on a set of indices with density $\geq 1 - \varepsilon$.

We always have $1 \leq \Gamma'_p(\{c_n\}) \leq \Gamma_p(\{c_n\}) \leq 2$. The sequence $c_n = n^\ell$ for $\ell \geq 2$ has $p$-adic collision exponent $2(1 - \ell^{-1})$, while its reduced collision exponent can be seen to be 1.

*Computation of $\Gamma'_q(\{c_n\})$ for a linear recursion.* We conclude this section with an algorithm for computing the reduced $q$-adic collision exponent of any linear recursion sequence and any integer $q > 1$. Let $\{c_k\}$ be such a sequence, i.e., for certain integers $a_0, \ldots, a_{n-1}$ ($a_0 \neq 0$) we have

(6) $\qquad$ for all $k > n$, $c_k + a_{n-1} c_{k-1} + a_{n-2} c_{k-2} + \cdots + a_0 c_{k-n} = 0$.

Denote by $f(x) = x^n + a_{n-1} x^{n-1} + \cdots + a_1 x + a_0$ the recursion polynomial of (6). We can assume that $f$ is of minimal degree. If $\{c_k\}$ is constant along some arithmetic progression, surely $\Gamma'_q(\{c_n\}) = 2$. Call $f$ *nondegenerate* if the only sequence $\{c_k\}$ satisfying $f$ and having a constant arithmetic subsequence is the zero sequence. We quote the following results:

THEOREM 4.3 ([13, Thms. 5.1–5.2]). *Let $f$ denote a linear recursion.*

(i) *$f$ is nondegenerate if and only if no roots of $f$, or their ratios, are roots of unity.*

(ii) *Let $\{c_k\}$ be a sequence of integers satisfying $f$, and suppose that $\{c_k\}$ has no constant arithmetic subsequences. Then $\Gamma'_q(\{c_k\}) = 1$ for any $q > 1$ relatively prime to $f(0)$.*

*Note.* In fact a stronger result is proved there: after discarding a set of arbitrarily small density, $\{c_k\}$ has *bounded cells*, i.e.,

$$\sup_n \max_{0 \leq t < q^n} |\{0 \leq k < q^n : c_k \equiv t \pmod{q^n}\}| < \infty.$$

For proofs of the results in the rest of this section, see Meiri [15].

Call a sequence $\{c_k\}$ an *intertwining* of the sequences $\{c_k^{(1)}\}, \ldots, \{c_k^{(N)}\}$ if $c_{N(k-1)+r} = c_k^{(r)}$ for all $k \geq 1$ and $r = 1, \ldots, N$.

LEMMA 4.4. *If $\{c_k\}$ is an intertwining of the sequences $\{c_k^{(1)}\}, \ldots, \{c_k^{(N)}\}$, then*

$$\Gamma_q(\{c_k\}) = \max_{1 \leq r \leq N} \Gamma_q(\{c_k^{(r)}\}), \tag{7}$$

*and a similar result holds for the reduced collision exponent.*

For a polynomial $f(x) = x^n + a_{n-1}x^{n-1} + \cdots + a_1 x + a_0$, define

$$\gcd[f] = \gcd(a_0, a_1, \ldots, a_{n-1}).$$

THEOREM 4.5. *Let $\{c_k\}$ be a linear recursion sequence, and suppose that its minimal recursion polynomial $f$ is nondegenerate. Decompose $q = q_1 q_2$ such that $q_1$ is the maximal factor of $q$ that is relatively prime to $\gcd[f]$. Then*

$$\Gamma'_q(\{c_k\}) = 2 - \frac{\log q_1}{\log q} = 1 + \frac{\log q_2}{\log q}. \tag{8}$$

*In particular we have $\Gamma'_q(\{c_k\}) = 1$ if and only if* no *prime factor of $q$ divides $\gcd[f]$, and $\Gamma'_q(\{c_k\}) = 2$ if and only if* all *prime factors of $q$ divide $\gcd[f]$.*

*Algorithm.* Computing the reduced $q$-adic collision exponent of a linear recursion sequence.

1. Compute the minimal recursion polynomial $f$ of $\{c_k\}$.
   *Method*: solve linear equations on the coefficients of $f$.

2. If $f(-1) = 0$, consider separately $\{c_{2k}\}$ and $\{c_{2k+1}\}$: each satisfies a recursion of lower degree; apply Lemma 4.4.

3. If $f(1) = 0, f'(1) \neq 0$, find constants $r, s$ such that the minimal polynomial of $\{rc_k - s\}$ does not vanish at 1, and replace $\{c_k\}$ by the latter sequence. The minimal polynomial of the new sequence does not vanish at 1.
   *Method*: find a rational number $s/r$ such that $\{c_k - s/r\}$ satisfies $f(x)/(x-1)$.

4. If $f(1) = f'(1) = 0$, then $\Gamma'_q(\{c_k\}) = 1$.

5. Check if $f$ is an intertwining of sequences, satisfying shorter recursions.
   *Method*: find the maximal number $D$ such that $\varphi(D) \leq n$, where $\varphi$ is Euler's totient function. For $d = 2, \ldots, D$, display $\{c_k\}$ as an intertwining of $d$ subsequences, and compute their minimal polynomials. If for some $d$ all the resulting polynomials have degree less than $\deg f$, apply Lemma 4.4. If not, then $f$ is nondegenerate.

6. If $f$ is nondegenerate, apply (8) to compute $\Gamma'_q(\{c_k\})$.

*Examples.*

1. The sequence $4, 9, 39, 219, \ldots$ which satisfies the recursion $c_k = 7c_{k-1} - 6c_{k-2}$ has minimal polynomial $f(x) = x^2 - 7x + 6$. In step 3 we see that $f(1) = 0$; indeed, $c_k = 3 + 6^k$ and so $\Gamma'_q(\{c_k\}) = \Gamma_q(\{c_k\}) = \Gamma_q(\{6^k\}) = 1 + \frac{\log 2^m 3^\ell}{\log q}$, assuming that $2^m \| q, 3^\ell \| q$. On the other hand, a sequence $\{c_k\}$ whose minimal polynomial is $(x-1)f(x)$ satisfies $\Gamma'_q(\{c_k\}) = 1$ for every $q > 1$.

2. Suppose that $f(x) = 1 + x^3 + x^6$ is the minimal recursion of $\{c_k\}$. Then $\{c_k\}$ is an intertwining of 3 sequences, each satisfying the Fibonacci recursion $F_k = F_{k-1} + F_{k-2}$. If either of these sequences is identically zero, then $\Gamma'_q(\{c_k\}) = 2$. Otherwise, $\Gamma'_q(\{c_k\}) = 1$.

3. If the minimal polynomial $f$ of $\{c_k\}$ is nondegenerate, and some prime factor of $q$ does not divide $\gcd[f]$, then $\Gamma'_q(\{c_k\}) < 2$.

4. Suppose that $c_k = r_1(k)q_1^k + \ldots + r_N(k)q_N^k$, where $r_i$ are polynomials. Suppose also that there exists some prime factor of $q$ that does not divide some $q_i$. Then $\Gamma'_q(\{c_k\}) < 2$.

## 5. Uniform distribution in subgroups

The next three simple lemmas are the key to proving Theorem 1.1.

LEMMA 5.1. *Let $\{X_n\}$ be an infinite sequence of independent random variables with values in $\mathbb{Z}_N \stackrel{\text{def}}{=} \mathbb{Z}/N\mathbb{Z}$, for some fixed integer $N > 1$. Suppose that for some nonzero $g \in \mathbb{Z}_N$,*

$$(9) \qquad \sum_{j=1}^{\infty} \sum_{x=0}^{N-1} \min\{\mathbb{P}(X_j = x), \mathbb{P}(X_j = x + g)\} = \infty.$$

*Let $S_n = X_1 + \cdots + X_n \pmod{N}$. Then for any $x \in \mathbb{Z}_N$,*

$$(10) \qquad \lim_{n \to \infty} \left( \mathbb{P}(S_n = x + g) - \mathbb{P}(S_n = x) \right) = 0.$$

*Proof.* We first prove the following property of Fourier coefficients of $S_n$:

$$(11) \qquad \text{for all } \ell \in \mathbb{Z}_N, \quad g\ell \not\equiv 0 \pmod{N} \quad \Longrightarrow \quad \lim_{n \to \infty} \widehat{S}_n(\ell) = 0.$$

For some $x \in \mathbb{Z}_N$,

$$\sum_{j=1}^{\infty} \min\{\mathbb{P}(X_j = x), \mathbb{P}(X_j = x + g)\} = \infty.$$

Set $p_j = \min\{\mathbb{P}(X_j = x), \mathbb{P}(X_j = x + g)\}$. Denoting $\varphi_\ell(t) = \exp(2\pi i \ell t/N)$, we have

$$\widehat{S}_n(\ell) = \mathbf{E}\,\varphi_\ell(\sum_{j=1}^{n} X_j) = \prod_{j=1}^{n} \mathbf{E}\,\varphi_\ell(X_j).$$

Write

$$\begin{aligned}
|\mathbf{E}\,\varphi_\ell(X_j)| &= \left|\sum_{k=0}^{N-1} \mathbb{P}(X_j = k)\varphi_\ell(k)\right| \\
&\leq p_j \left|\frac{\varphi_\ell(x) + \varphi_\ell(x+g)}{2}\right| + (1 - p_j) \\
&= p_j \left|\frac{1 + \exp(2\pi i \ell g/N)}{2}\right| + (1 - p_j).
\end{aligned}$$

Since $N \nmid g\ell$, we have $|1+\exp(2\pi i \ell g/N)|/2 \leq \gamma < 1$, for $\gamma \stackrel{\text{def}}{=} |1+\exp(2\pi i/N)|/2$. Hence

$$|\mathbf{E}\,\varphi_\ell(X_j)| \leq 1 - (1-\gamma)p_j.$$

By our assumptions, $\sum_{j=1}^{\infty}(1-\gamma)p_j = \infty$; hence $\lim_{n\to\infty} \prod_{j=1}^{n} |\mathbf{E}\,\varphi_\ell(X_j)| = 0$, and (11) follows.

To prove (10), use inverse Fourier transform to write

$$\mathbb{P}(S_n = x + g) - \mathbb{P}(S_n = x) = \frac{1}{N} \sum_{\ell=0}^{N-1} \widehat{S}_n(\ell)\Big(\varphi_{-\ell}(x+g) - \varphi_{-\ell}(x)\Big).$$

If $N | g\ell$ we have $\varphi_{-\ell}(x+g) - \varphi_{-\ell}(x) = 0$. If $N \nmid g\ell$, apply (11). $\square$

LEMMA 5.2. *Let $\{X_n\}$ be an infinite sequence of independent random variables with values in $\mathbb{Z}_N \stackrel{\text{def}}{=} \mathbb{Z}/N\mathbb{Z}$, for some fixed integer $N > 1$. Suppose that there exists a subgroup $G \subseteq \mathbb{Z}_N$, generated by $g_1, \ldots, g_r$, such that (9) holds for $g = g_1, \ldots, g_r$. Let $S_n = X_1 + \cdots + X_n \pmod{N}$, and let $S_n \bmod G$ denote the projection of $S_n$ to $\mathbb{Z}_N/G$. Then*

$$\mathbf{E}\,H(S_n | S_n \bmod G) \longrightarrow \log|G|.$$

*Proof.* Let $g_0 \in G$ be a generator of $G$. Applying Lemma 5.1 for $g = g_1, \ldots, g_r$, we see that (10) holds for $g = g_0$ as well. Denote for a moment $\widetilde{S}_n = S_n \bmod G$. If $G = \mathbb{Z}_N$, the result is clear: $\widetilde{S}_n$ is constant, while the distribution of $S_n$ converges to the uniform distribution on $\mathbb{Z}_N$, whose entropy is $\log|\mathbb{Z}_N| = \log|G|$. In the general case, write

$$\mathbf{E}\,H(S_n|\widetilde{S}_n) = \sum_{t+G\in \mathbb{Z}_N/G} \mathbb{P}(\widetilde{S}_n = t+G) H(S_n|\widetilde{S}_n = t+G).$$

By (10) for $g = g_0$, the conditional distribution of $S_n$ on every coset $t+G$ converges to the uniform distribution whose entropy is again $\log|G|$, and we are done. □

In the next lemma we restrict our attention to $N = p^k$ and identify $\mathbb{Z}_{p^k}$ with $k$-digit numbers in base $p$. For $n \leq k$, let $\pi_n : \mathbb{Z}_{p^k} \longrightarrow \mathbb{Z}_{p^n}$ denote the projection to the $n$ most-significant digits, i.e., $\pi_n(z) = \lfloor z/p^{k-n} \rfloor$.

LEMMA 5.3. *Let $Y$ be a $\mathbb{Z}_{p^k}$-valued random variable. Let $G$ be a subgroup of $\mathbb{Z}_{p^k}$, and suppose that $n$ satisfies $p^n \geq |G|$. Then*

$$H(\pi_n(Y)) \geq H(Y|Y \bmod G) = H(Y) - H(Y \bmod G).$$

*Proof.* In every coset $t+G$, the projection $\pi_n$ is one-to-one. Thus $\pi_n(Y)$ and $Y \bmod G$ determine $Y$; hence $H(\pi_n(Y)) + H(Y \bmod G) \geq H(Y)$. □

Finally, we state the following simple lemma that will be used to show that the subgroup $G$ we get in 5.1 will be large. For $t = \sum_{i=1}^{k} t_i p^{k-i} \in \mathbb{Z}/p^k\mathbb{Z}$ we call $t_i$ the $i^{\text{th}}$ digit of $t$.

LEMMA 5.4. *Let $t \in G \subset \mathbb{Z}/p^k\mathbb{Z}$, and suppose that the $i^{\text{th}}$ digit of $t$ is nonzero. Then $t$ generates in $G$ a subgroup of order at least $p_*^i$, where $p_*$ is the smallest prime divisor of $p$.*

We leave the proof of this lemma to the reader.

## 6. Entropy and subgroups

*Notation.* Suppose that $\tilde{\mu}$ is a $\sigma_p \times \sigma_p \times \cdots$-invariant measure on $\mathbb{T}^{\mathbb{N}}$. For $x \in \mathbb{T}^{\mathbb{N}}$ we denote its projection on coordinates by superscripts, e.g., $x^{a\ldots b} = \pi^{a\ldots b}(x)$ denotes projection to coordinates $a,\ldots,b$. Similarly, $x^a = x^{a\ldots a} = \pi^a(x)$. We identify numbers in $\mathbb{T}$ with their base-$p$ expansion. In this identification, $\sigma_p$ is equivalent to a shift map on $\{0,\ldots,p-1\}^{\mathbb{N}}$, and we denote by $x^i_{c\ldots d}$ the $d-c+1$ digits starting at point $c$ in the base-$p$ expansion of $x^i$. For every $i$ let $\alpha_1^i$ denote the partition of $\mathbb{T}^{\mathbb{N}}$ according to the first digit of $x^i$, and we use the notation

$$\alpha^{a\ldots b}_{c\ldots d} \stackrel{\text{def}}{=} \bigvee_{i=a}^{b} \bigvee_{j=c}^{d} \sigma_p^{-j+1} \alpha_1^i$$

for the partition of $\mathbb{T}^{\mathbb{N}}$ according to the digits $c, \ldots, d$ of $x^a, \ldots, x^b$. By a slight abuse of notation $\alpha_1^i$ also denotes the analogous partition of $\pi^i(\mathbb{T}^{\mathbb{N}})$. Recall the definition of the function $\Theta^n \colon \mathbb{T}^{\mathbb{N}} \to \mathbb{T}$ as the sum mod 1 of $x^1, \ldots, x^n$. We denote the partition of $\mathbb{T}^{\mathbb{N}}$ according to the first $k$ digits of $\Theta^n(x)$ by $\theta_{1\ldots k}^n$, etc.

We also use the following notation for a measure $\mu$ and a finite partition $\gamma$: we denote by $D_\mu(\gamma)$ the probability vector $\{\mu(C) : C \in \gamma\}$, and by $H_\mu(\gamma) = H(D_\mu(\gamma))$ its entropy. Similarly, given a $\sigma$-algebra $\mathfrak{B}$ we denote by $D_\mu(\gamma|\mathfrak{B})$ the conditional distribution of $\gamma$ given $\mathfrak{B}$, and by $H_\mu(\gamma|\mathfrak{B}) = H(D_\mu(\gamma|\mathfrak{B}))$ its entropy; note that this is a function on the measure space, not a constant.

In our context it is natural to view $x_{1\ldots k}^a$ as an element of $\mathbb{Z}/p^k\mathbb{Z}$, and thus if $G$ is a subgroup of $\mathbb{Z}/p^k\mathbb{Z}$, it makes sense to consider $x_{1\ldots k}^a \bmod G$ or the corresponding partition of $\pi^a(\mathbb{T}^{\mathbb{N}})$, which is $\alpha_{1\ldots k}^a \bmod G$.

PROPOSITION 6.1. *Let $\mu, \nu$ be two $p$-invariant measures, and $G$ a subgroup of $\mathbb{Z}/p^k\mathbb{Z}$ for some $k \in \mathbb{N}$. Then*

$$\mathbf{E}\, H_{\mu * \nu}(x_{1\ldots k} \bmod G | x_{k+1\ldots\infty}) \geq \mathbf{E}\, H_\mu(x_{1\ldots k} \bmod G | x_{k+1\ldots\infty}).$$

*Proof.* Define independent random variables $X, Y$ with $X \sim \mu$ and $Y \sim \nu$. Then $X + Y \sim \mu * \nu$. Denote by $\alpha, \beta, \gamma$ the partitions according to the first digit in the base-$p$ expansion of $X, Y$ and $X + Y$, respectively. We have

(12)
$$\begin{aligned} \mathbf{E}\, H_{\mu * \nu}(x_{1\ldots k} \bmod G | x_{k+1\ldots\infty}) &= \mathbf{E}\, H_{\mu \times \nu}(\gamma_{1\ldots k} \bmod G | \gamma_{k+1\ldots\infty}) \\ &\geq \mathbf{E}\, H_{\mu \times \nu}(\gamma_{1\ldots k} \bmod G | \alpha_{k+1\ldots\infty} \vee \beta_{1\ldots\infty}). \end{aligned}$$

However, given all we have conditioned upon on the right-hand side of (12), $\alpha_{1\ldots k} \bmod G$ uniquely determines $\gamma_{1\ldots k} \bmod G$, and vice versa. Hence

$$\begin{aligned} H_{\mu \times \nu}(\gamma_{1\ldots k} \bmod G | \alpha_{k+1\ldots\infty} \vee \beta_{1\ldots\infty}) &= H_{\mu \times \nu}(\alpha_{1\ldots k} \bmod G | \alpha_{k+1\ldots\infty} \vee \beta_{1\ldots\infty}) \\ &= H_{\mu \times \nu}(\alpha_{1\ldots k} \bmod G | \alpha_{k+1\ldots\infty}), \end{aligned}$$

using the independence of $\alpha$ and $\beta$ in the last equality. This gives the required inequality. $\square$

COROLLARY 6.2. *Let $\{\mu_i\}$ be a sequence of $p$-invariant measures, and denote $\tilde{\mu} = \prod \mu_i$. Suppose that $G$ is a subgroup of $\mathbb{Z}/p^k\mathbb{Z}$ for some $k \in \mathbb{N}$. Then*

$$\mathbf{E}\, H_{\tilde{\mu}}(\theta_{1\ldots k}^n \bmod G | \theta_{k+1\ldots\infty}^n)$$

*is monotone nondecreasing in $n$.*

LEMMA 6.3. *Let $\mu$ be a measure on $\mathbb{T}$, and suppose that $G \subset \mathbb{Z}/p^k\mathbb{Z}$ is a group of size $\geq p^\ell$. Then*

$$H(\alpha_{1\ldots\ell}) \geq (\ell - 1)\log p - \log|G| + \int H(\alpha_{1\ldots k} | \alpha_{1\ldots k} \bmod G \vee \alpha_{k+1\ldots\infty})\, d\mu.$$

*Proof.* Let $n = \lceil \log |G| / \log p \rceil$. Applying Lemma 5.3 in each fiber we get

$$H(\alpha_{1...n}) \geq \int H(\alpha_{1...n}|\alpha_{k+1...\infty})\,d\mu \geq \int H(\alpha_{1...k}|\alpha_{1...k} \bmod G \vee \alpha_{k+1...\infty})\,d\mu.$$

Since $n \geq \ell$, we also have

$$n \log p - H(\alpha_{1...n}) \geq \ell \log p - H(\alpha_{1...\ell}).$$

Now combine the two inequalities. $\square$

Finally, we quote the following extension of the Borel-Cantelli Lemma ([2, ex. 22.4]).

LEMMA 6.4. *Let $Y_1, Y_2, \ldots$ be a sequence of nonnegative bounded independent random variables. Suppose that $\sum \mathbf{E} Y_i = \infty$. Then $\sum Y_i = \infty$ almost-surely.*

## 7. The Convolution Theorem

In this section we prove Theorem 1.1. Fix $p$, and consider the function

$$(13) \qquad \psi(\beta) = H(1-\beta, \frac{\beta}{p-1}, \ldots, \frac{\beta}{p-1}).$$

Clearly $\psi : [0, 1 - \frac{1}{p}] \to [0, \log p]$ is increasing, onto and concave. A quick calculation shows that the inverse function $\psi^{-1}$ satisfies

$$(14) \qquad \psi^{-1}(h) \geq C \cdot \frac{h}{|\log(h/\log p)|}$$

for some constant $C$ and for all $h$ sufficiently small.

Given a probability distribution $\nu$ on $\{0, \ldots, p-1\}$ define

$$\|\nu\|_\infty = \max_{k=0,\ldots,p-1} \nu(k).$$

It can be easily verified that

$$(15) \qquad 1 - \|\nu\|_\infty \leq \frac{1}{p} \sum_{x=0}^{p-1} \sum_{g=1}^{p-1} \min\{\nu(x), \nu(x+g \bmod p)\}$$

and also

$$(16) \qquad H(\nu) \leq \psi(1 - \|\nu\|_\infty).$$

Let $\mu$ denote an nonatomic measure on $\mathbb{T}$, not necessarily invariant. For $x \in \mathbb{T}$ denote by $x = \sum x_i p^{-i}$ its expansion in base $p$, and define

$$I_k(\mu)(x) = -\log \mu(x_{1...k}|x_{k+1...\infty}).$$

Then $\mathbf{E}\, I_k(\mu) = \mathbf{E}\, H_\mu(x_{1\ldots k}|x_{k+1\ldots\infty})$. Similarly, given a subgroup $G \subset \mathbb{Z}_{p^k} \stackrel{\text{def}}{=} \mathbb{Z}/p^k\mathbb{Z}$ define

$$I_G(\mu)(x) = -\log \mu(x_{1\ldots k}|x_{1\ldots k} \bmod G \vee x_{k+1\ldots\infty}).$$

Then $I_k = I_G$ for $G = \mathbb{Z}_{p^k}$. We remark that

$$\mathbf{E}\, I_G(\mu) = -\int \log\left(\frac{d\mu}{d\delta_G * \mu}\right) d\mu,$$

where $\delta_G \stackrel{\text{def}}{=} \sum_{g \in G} \delta_{g/p^k}$.

LEMMA 7.1. *Let $G \subset \mathbb{Z}_{p^k}$, and suppose that $\mu$ and $\nu$ are nonatomic measures on $\mathbb{T}$. Then*

$$\mathbf{E}\, I_k(\mu * \nu) - \mathbf{E}\, I_k(\mu) \geq \mathbf{E}\, I_G(\mu * \nu) - \mathbf{E}\, I_G(\mu).$$

*Proof.* Rewrite the desired inequality as

(17) $$\mathbf{E}\, I_k(\mu * \nu) - \mathbf{E}\, I_G(\mu * \nu) \geq \mathbf{E}\, I_k(\mu) - \mathbf{E}\, I_G(\mu).$$

We have

$$\begin{aligned}
\mathbf{E}\, I_k(\mu) - \mathbf{E}\, I_G(\mu) &= \mathbf{E}\, H_\mu(x_{1\ldots k}|x_{k+1\ldots\infty}) \\
&\quad - \mathbf{E}\, H_\mu(x_{1\ldots k}|x_{1\ldots k} \bmod G \vee x_{k+1\ldots\infty}) \\
&= \mathbf{E}\, H_\mu(x_{1\ldots k} \bmod G|x_{k+1\ldots\infty}).
\end{aligned}$$

Applying the same consideration to $\mu * \nu$, we see that (17) is equivalent to

$$\mathbf{E}\, H_{\mu*\nu}(x_{1\ldots k} \bmod G|x_{k+1\ldots\infty}) \geq \mathbf{E}\, H_\mu(x_{1\ldots k} \bmod G|x_{k+1\ldots\infty}),$$

which follows from Proposition 6.1. $\square$

LEMMA 7.2. *Let $\{\mu_n\}$ be a sequence of probability measures on $\mathbb{T}$, and form the product measure $\tilde{\mu} = \prod \mu_i$. Suppose that for some fixed number $k$,*

(18) $$\sum_{i=1}^{\infty} \mathbf{E}\, \psi^{-1}\left(H_{\mu_i}(x_k^i|t_{k+1\ldots\infty}^i)\right) = \infty.$$

*Then for $\tilde{\mu}$-almost-every $t \in \mathbb{T}^{\mathbb{N}}$ there exists a group $G_k(t) \subset \mathbb{Z}_{p^k}$ such that*

(19) $$\mathbf{E}\, I_{G_k(t)}(\mu_1 * \cdots * \mu_n) \longrightarrow \mathbf{E}\, \log|G_k(t)|.$$

*Furthermore, the map $t \mapsto G_k(t)$ is measurable, and $|G_k(t)| \geq p_*^k$ a.e., where $p_*$ is the smallest prime factor of $p$.*

*Proof.* We claim that

(20) $$\sum_{n=1}^{\infty} \sum_{x_k=0}^{p-1} \sum_{g_k=1}^{p-1} \min\{\mu_n(x_k|t_{k+1\ldots\infty}^n), \mu_n(x_k + g_k \bmod p|t_{k+1\ldots\infty}^n)\} = \infty \quad \tilde{\mu}\text{-a.e.}$$

Indeed, by Lemma 6.4, for every $t$ in a set of full $\tilde{\mu}$ measure,

$$\sum_{i=1}^{\infty} \psi^{-1}\left(H_{\mu_i}(x_k^i | t_{k+1\ldots\infty}^i)\right) = \infty.$$

By (16), for all such $t$ we have

$$\sum_{n=1}^{\infty} \left(1 - \|D_{\mu_n}(\alpha_k^n | t_{k+1\ldots\infty}^n)\|_\infty\right) = \infty.$$

Using (15), we conclude that a.e. $t$ satisfies (20).

Define for $t \in \mathbb{T}^{\mathbb{N}}$ a group $G_k(t)$ by

$$G_k(t) = \left\langle g \in \mathbb{Z}_{p^k} \;\middle|\; \sum_{n=1}^{\infty} \sum_{x=0}^{p^k-1} \min\{\mu_n(x|t_{k+1\ldots\infty}^n), \mu_n(x+g \bmod p^k | t_{k+1\ldots\infty}^n)\} = \infty \right\rangle.$$

To prove the lemma, take some $t$ satisfying (20). We first show that $|G_k(t)| \geq p_*^k$. There exist distinct $x_k, y_k \in \{0, \ldots, p-1\}$ such that

$$\sum_{n=1}^{\infty} \mu_n(x_k | t_{k+1\ldots\infty}^n) = \infty,$$

and similarly for $y_k$. Since

$$\mu_n(x_k | t_{k+1\ldots\infty}^n) = \sum_{\substack{x=0,\ldots,p^k-1 \\ x \equiv x_k \;(\bmod\; p)}} \mu_n(x | t_{k+1\ldots\infty}^n),$$

there exists $x \in \{0, \ldots, p^k - 1\}$ with $x \equiv x_k \pmod{p}$ such that

$$\sum_{n=1}^{\infty} \mu_n(x | t_{k+1\ldots\infty}^n) = \infty.$$

Similarly, there exists $y \in \{0, \ldots, p^k - 1\}$ with $y \equiv y_k \pmod{p}$ and

$$\sum_{n=1}^{\infty} \mu_n(y | t_{k+1\ldots\infty}^n) = \infty.$$

By definition, $y - x \in G_k(t)$. But the least significant digit of $y - x$ is $y_k - x_k \neq 0$. By Lemma 5.4, $|G_k(t)| \geq p_*^k$.

It remains to prove (19). Clearly, $I_{G_k(t)}(\mu_1 * \cdots * \mu_n) \leq \log |G_k(t)|$. But

$$\begin{aligned}
\mathbf{E}\, I_{G_k(t)}(\mu_1 * \cdots * \mu_n) &= \mathbf{E}\, H_{\tilde{\mu}}(\theta_{1\ldots k}^n | \theta_{1\ldots k}^n \bmod G_k(t) \vee \theta_{k+1\ldots\infty}^n) \\
&\geq \mathbf{E}\, H_{\tilde{\mu}}(\theta_{1\ldots k}^n | \theta_{1\ldots k}^n \bmod G_k(t) \vee t_{k+1\ldots\infty}^{1\ldots n}) \\
&\longrightarrow \mathbf{E} \log |G_k(t)|,
\end{aligned}$$

by an application of Lemma 5.2 for the random variables $X_i \sim D_{\mu_i}(\alpha_{1\ldots k}^i | t_{k+1\ldots\infty}^i)$. This concludes the proof of (19). □

LEMMA 7.3. *Under the assumptions of Lemma 7.2, define*

$$h_k = \sup_n \frac{1}{k} \mathbf{E}\, H_{\mu_1 * \cdots * \mu_n}(x_{1\ldots k}|x_{k+1\ldots\infty}). \tag{21}$$

For any $m$, if $\frac{1}{k}\mathbf{E}\, H_{\mu_1 * \cdots * \mu_m}(x_{1\ldots k}|x_{k+1\ldots\infty}) > h_k - \varepsilon$, then

$$H_{\mu_1 * \cdots * \mu_m}(t_{1\ldots \ell}) \geq (\ell-1)\log p - (k+1)\varepsilon, \tag{22}$$

where $\ell = \ell(k) = \lfloor k \frac{\log p_*}{\log p} \rfloor$.

*Proof.* We have $|G_k(t)| \geq p_*^k \geq p^\ell$ almost-everywhere. By Lemma 6.3,

$$H_{\mu_1 * \cdots * \mu_m}(t_{1\ldots\ell}) \geq (\ell-1)\log p - \int \left[ \log|G_k| - I_{G_k}(\mu_1 * \cdots * \mu_m) \right] d\tilde{\mu}. \tag{23}$$

By Lemma 7.2, there exists some $n > m$ such that $\mathbf{E}\, I_{G_k}(\mu_1 * \cdots * \mu_n) \geq \mathbf{E}\log|G_k| - \varepsilon$. Applying Lemma 7.1 we get

$$\begin{aligned}
k\varepsilon &\geq \mathbf{E}\, I_k(\mu_1 * \cdots * \mu_n) - \mathbf{E}\, I_k(\mu_1 * \cdots * \mu_m) \\
&\geq \mathbf{E}\, I_{G_k}(\mu_1 * \cdots * \mu_n) - \mathbf{E}\, I_{G_k}(\mu_1 * \cdots * \mu_m) \\
&\geq \mathbf{E}\log|G_k| - \varepsilon - \mathbf{E}\, I_{G_k}(\mu_1 * \cdots * \mu_m).
\end{aligned}$$

By substituting this into (23) we obtain the desired result. □

*Proof of Theorem* 1.1. By convexity of $\psi^{-1}$,

$$\begin{aligned}
\sum_{n=1}^{\infty} \mathbf{E}\, \psi^{-1}\left(H_{\mu_n}(t_k|t_{k+1\ldots\infty})\right) &\geq \sum_{n=1}^{\infty} \psi^{-1}\left(\mathbf{E}\, H_{\mu_n}(t_k|t_{k+1\ldots\infty})\right) \\
&= \sum_{n=1}^{\infty} \psi^{-1}\left(h(\mu_n, \sigma_p)\right) = \infty,
\end{aligned}$$

applying (1) and (14). Hence the assumptions of Lemmas 7.2–7.3 hold.

By Proposition 6.2, $h(\mu_1 * \cdots * \mu_n, \sigma_p)$ is monotone nondecreasing in $n$. Define

$$h \stackrel{\text{def}}{=} \lim_{n\to\infty} h(\mu_1 * \cdots * \mu_n, \sigma_p) = \sup_n h(\mu_1 * \cdots * \mu_n, \sigma_p).$$

We need to show that $h = \log p$.

The $h_k$ defined in (21) satisfy $h_k = h$ for all $k$, by invariance of $\mu_1 * \cdots * \mu_n$. For arbitrary $\varepsilon$, let $m$ be big enough so that $h(\mu_1 * \cdots * \mu_m, \sigma_p) > h - \varepsilon$. By Lemma 7.3, (22) holds for every $k$, with $\ell(k) = \lfloor k\frac{\log p_*}{\log p}\rfloor$. Hence

$$h(\mu_1 * \cdots * \mu_m, \sigma_p) \geq \lim_{k\to\infty} \frac{(\ell(k)-1)\log p - (k+1)\varepsilon}{\ell(k)} = \log p - \varepsilon \frac{\log p}{\log p_*}.$$

Letting $\varepsilon \downarrow 0$ proves the theorem. □

## 8. Joinings of full entropy

In this section we study basic properties of joinings of full entropy, and prove Theorem 1.8.

Let $\alpha^{a...b} \stackrel{\text{def}}{=} \alpha^{a...b}_{1...\infty}$, and denote by $\mathcal{T}^i \stackrel{\text{def}}{=} \bigwedge_{j=1}^{\infty} \alpha^i_{j...\infty}$ the tail $\sigma$-algebra of $\pi^i(\mathbb{T}^{\mathbb{N}})$. We denote by $\mathcal{T}$ the $\sigma$-algebra $\bigvee_{a=1}^{\infty} \left( \bigwedge_{j=1}^{\infty} \alpha^{1...a}_{j...\infty} \right)$.

The simplest example of a joining of full entropy is a product measure $\tilde{\mu} = \prod \mu_i$. The following lemma shows that in the general case, a similar independence property holds.

LEMMA 8.1. *Let $\tilde{\mu}$ be a joining of full entropy of $\{\mu_i\}_{i=1}^{\infty}$. Then the random variables $\{x^i\}$ are conditionally independent given $\mathcal{T}$.*

*Proof.* For any $k, n \geq 1$ we have

$$(24) \quad \int H_{\tilde{\mu}}(\alpha^{1...n}_{1...k} | \alpha^{1...n}_{k+1...\infty}) \, d\tilde{\mu} \leq \int \sum_{i=1}^{n} H_{\tilde{\mu}}(\alpha^i_{1...k} | \alpha^{1...n}_{k+1...\infty}) \, d\tilde{\mu}$$

$$(25) \quad \leq \int \sum_{i=1}^{n} H_{\mu_i}(\alpha^i_{1...k} | \alpha^i_{k+1...\infty}) \, d\mu_i.$$

Since $h(\pi^{1...n}(\tilde{\mu}), \sigma_p \times \cdots \times \sigma_p) = \sum_{i=1}^{n} h(\mu_i, \sigma_p)$, both inequalities must in fact be equalities. Equality in (24) implies that the random variables $\{x^i_{1...k}\}_{i=1}^{n}$ are independent given $\alpha^{1...n}_{k+1...\infty}$. A reverse martingale argument shows that in fact $\{x^i_{1...k}\}_{i=1}^{n}$ are independent given $\bigwedge_{j=1}^{\infty} \alpha^{1...n}_{j...\infty}$. By a standard martingale argument, $\{x^i_{1...k}\}_{i=1}^{n}$ are independent given $\mathcal{T} = \bigvee_{a=1}^{\infty} \left( \bigwedge_{j=1}^{\infty} \alpha^{1...a}_{j...\infty} \right)$. Since this is true for any $k, n \geq 1$, the assertion follows. □

We obtain the following corollary:

LEMMA 8.2. *Let $\tilde{\mu}$ denote a joining of full entropy of $\{\mu_i\}_{i=1}^{\infty}$, and denote by $\mathcal{L} = \bigwedge_{a=1}^{\infty} \alpha^{a...\infty}_{1...\infty}$ the tail $\sigma$-algebra of the stochastic process $\{x^i\}$. Then $\mathcal{L} \subset \mathcal{T}$. ("The bottom tail is contained in the right tail.")*

*Proof.* Let $A$ be some $\mathcal{L}$-measurable set. Denote by $\tau$ the factor map corresponding to $\mathcal{T}$. By Lemma 8.1, given $\tau(x)$ the random variables $\{x^i\}_{i=1}^{\infty}$ are independent. By Kolmogorov's 0-1 Law, $\tilde{\mu}(A|\mathcal{T})$ is either 0 or 1. □

We cite the following technical lemma:

LEMMA 8.3. *Let $(X, T, \mathcal{D}, \nu)$ denote a measure preserving system, and suppose that $\beta$ and $\gamma$ are two finite partitions of $X$. Denote by $\mathcal{T}_\gamma$ the tail*

$\sigma$-algebra corresponding to $\gamma$ and $\nu$. Then for all $k \geq 1$,

$$D(\beta_{1...k}|\beta_{k+1...\infty} \vee \mathcal{T}_\gamma) = D(\beta_{1...k}|\beta_{k+1...\infty}) \quad \nu-\text{a.e.}$$

The same applies if instead of $\mathcal{T}_\gamma$ one takes a limit of an increasing sequence of tails of finite partitions.

*Proof.* See [17, Cor. 5.28, p. 99], or [16, Lemma 7, p. 65]. □

We also need the following monotonicity property of conditional entropy, which follows from Jensen's inequality.

LEMMA 8.4. *Let $\alpha, \beta, \beta', \gamma$ denote partitions in a measure preserving system $(X, T, \mathcal{D}, \nu)$, and suppose that $\gamma \leq \beta \leq \beta'$. Then*

$$\mathbf{E}\left(H_\nu(\alpha|\beta) \,\bigg|\, \gamma\right)(x) \geq \mathbf{E}\left(H_\nu(\alpha|\beta') \,\bigg|\, \gamma\right)(x) \quad \nu-\text{a.e.}$$

*The same applies if $\beta$ and $\gamma$ are limits of increasing sequences of partitions.*

PROPOSITION 8.5. *Let $\tilde{\mu}$ be a joining of full entropy, $k \in \mathbb{N}$, and $G$ a subgroup of $\mathbb{Z}/p^k\mathbb{Z}$. Then for a.e. $x$,*

$$\mathbf{E}\left(H_{\tilde{\mu}}(\theta_{1...k}^n \bmod G|\theta_{k+1...\infty}^n) \,\bigg|\, \mathcal{T}\right)(x)$$

*is monotone nondecreasing in $n$.*

*Proof.* By Lemma 8.3, we have $\tilde{\mu}$-a.e.

$$\mathbf{E}\left(H_{\tilde{\mu}}(\theta_{1...k}^n \bmod G|\theta_{k+1...\infty}^n) \,\bigg|\, \mathcal{T}\right)(x)$$
$$= \mathbf{E}\left(H_{\tilde{\mu}}(\theta_{1...k}^n \bmod G|\theta_{k+1...\infty}^n \vee \mathcal{T}) \,\bigg|\, \mathcal{T}\right)(x).$$

Now apply 6.2 in each fiber. □

We leave the proof of the following estimate to the reader.

LEMMA 8.6. *Let $X$ be a random variable such that $0 \leq X \leq M$ and $\mathbf{E} X \geq \eta$. Then*

$$\mathbb{P}(X > \frac{\eta}{2}) \geq \frac{\eta}{2M}.$$

COROLLARY 8.7. *If $X$ is as in Lemma 8.6, $g$ a monotone nondecreasing function, then*

$$\mathbf{E}(g(X)) \geq \frac{\eta g(\eta/2)}{2M}.$$

The last lemma we need is the following variant of Shannon's theorem. Denote by $(\mu|F)$ the measure $\mu$ conditioned on $F$, i.e.,

$$(\mu|F)(A) = \mu(A \cap F)/\mu(F).$$

LEMMA 8.8. *Let $\beta$ be a finite partition in an ergodic measure preserving system $(X, \mu, T)$, and denote $\beta_{1\ldots k} = \beta \vee \ldots \vee T^{-(k-1)}\beta$. Suppose that for some sequence of subsets $\{F_k\}$ and a constant $\gamma > 0$, for every $k$ large enough $\mu(F_k) \geq \gamma$. Then*

$$h_\mu(\beta, T) \geq \limsup_{k \to \infty} \frac{1}{k} H_{\mu|F_k}(\beta_{1\ldots k}).$$

*Proof.* Let $\varepsilon > 0$, and denote by $N$ the number of atoms in $\beta$. By Shannon's theorem, for $k$ large enough it is possible to cover a subset of $X$ of $\mu$-measure $1 - \varepsilon\gamma$ by at most $\exp\left(k(h_\mu(\beta, T) + \varepsilon)\right)$ atoms of $\beta_{1\ldots k}$. In particular, the $\mu|F_k$-measure of the union of these atoms is at least $1 - \varepsilon$. Hence by a standard property of entropy (see Rudolph [17, Cor. 5.17])

$$\frac{1}{k} H_{\mu|F_k}(\beta_{1\ldots k}) \leq \frac{1}{k} \cdot \left(H(\varepsilon, 1 - \varepsilon) + k(h_\mu(\beta, T) + \varepsilon) + \varepsilon k \log N\right).$$

Letting $k \to \infty$ and $\varepsilon \to 0$ proves the lemma. □

*Proof of Theorem* 1.8. Arguing as in Lemma 2.1, the ergodic components of $\tilde{\mu}$ are also joinings of full entropy of $\{\mu_i\}$. Since by Rokhlin's theorem the entropy of a measure is the average of the entropies of its ergodic components, it is enough to consider the case where $\tilde{\mu}$ is ergodic.

By Proposition 8.5, $h(\Theta^n \tilde{\mu}, \sigma_p)$ is monotone nondecreasing in $n$. Define

$$h \stackrel{\mathrm{def}}{=} \lim_{n \to \infty} h(\Theta^n \tilde{\mu}, \sigma_p) = \sup_n h(\Theta^n \tilde{\mu}, \sigma_p).$$

We need to show that $h = \log p$.

Define conditional measures $\bar{\mu}_i \sim D_{\tilde{\mu}}(\alpha^i_{1\ldots\infty}|\mathcal{T})$. The sequence $\{\bar{\mu}_i\}$ is a random sequence of measures, depending on a choice of fiber (a point in $\mathcal{T}$). By Lemma 8.1, $\{\bar{\mu}_i\}$ are independent. We claim that the assumptions of Lemmas 7.2–7.3 hold for $\{\bar{\mu}_i\}$ with positive probability. To see this, consider the random variables

$$Z_i = \mathbf{E}\left(\psi^{-1}\left(H_{\mu_i}(\alpha^i_k | t^i_{k+1\ldots\infty})\right) \Big| \mathcal{T}\right)$$

$$W_i = \mathbf{E}\left(H_{\mu_i}(\alpha^i_k | t^i_{k+1\ldots\infty}) \Big| \mathcal{T}\right).$$

To prove that (18) holds with positive probability, we need to show that $\mathbb{P}(\sum Z_i = \infty) > 0$. Let $\eta \stackrel{\mathrm{def}}{=} \inf_i h(\mu_i, \sigma_p) > 0$. By Corollary 8.7

(26) $$Z_i \geq \frac{W_i \psi^{-1}(W_i/2)}{2 \log p}.$$

Since $\mathbf{E}(W_i) = h(\mu_i, \sigma_p) \geq \eta$, by Lemma 8.6 and (26),

$$\mathbb{P}\left(Z_i > \frac{\eta \psi^{-1}(\eta/4)}{4 \log p}\right) > \eta^* \stackrel{\text{def}}{=} \frac{\eta}{2 \log p}.$$

Hence $\mathbb{P}(\sum Z_i = \infty) \geq \mathbb{P}(Z_i > \frac{\eta \psi^{-1}(\eta/4)}{4 \log p} \text{ i.o.}) \geq \eta^* > 0$.

Fix $\varepsilon \in (0, \eta^*)$, and choose $m$ with $h - h(\Theta^m \tilde{\mu}, \sigma_p) \leq \varepsilon^2$. For any $k$, define

$$\Upsilon_k^{(n)}(t) = \frac{1}{k} \mathbf{E}\left(H_{\tilde{\mu}}(\theta_{1...k}^n | \theta_{k+1...\infty}^n) \,\bigg|\, \mathcal{T}\right)(t),$$

$$\Upsilon_k(t) = \sup_n \Upsilon_k^{(n)}(t).$$

By Proposition 8.5, $\Upsilon_k^{(n)}(t)$ is almost-surely monotonically increasing in $n$. Thus

$$0 \leq \mathbf{E}(\Upsilon_k - \Upsilon_k^{(m)}) = h - h(\Theta^m \tilde{\mu}, \sigma_p) \leq \varepsilon^2;$$

hence by Markov's inequality

$$\tilde{\mu}[\Upsilon_k - \Upsilon_k^{(m)} \geq \varepsilon] \leq \frac{\varepsilon^2}{\varepsilon} = \varepsilon.$$

Define $F_k = \{\sum Z_i = \infty\} \setminus \{\Upsilon_k - \Upsilon_k^{(m)} \geq \varepsilon\}$. The event $\{\sum Z_i = \infty\}$ is $\mathcal{T}$-measurable by Lemma 8.2; hence $F_k$ is $\mathcal{T}$-measurable as well. Clearly $\tilde{\mu}(F_k) \geq \eta^* - \varepsilon$. Applying Lemma 7.3, we have (denoting $\ell = \ell(k) = \lfloor k \frac{\log p_*}{\log p} \rfloor$)

$$\text{for all } t \in F_k, \quad H_{\tilde{\mu}}(\theta_{1...\ell}^m | \mathcal{T})(t) \geq (\ell - 1) \log p - (k+1)\varepsilon.$$

Thus by Lemma 8.4

$$H_{\tilde{\mu}|F_k}(\theta_{1...\ell}^m) = \mathbf{E}\left(H_{\tilde{\mu}}\left(\theta_{1...\ell}^m | \{F_k, F_k^c\}\right) \,\bigg|\, F_k\right) \geq (\ell - 1) \log p - (k+1)\varepsilon.$$

Dividing by $\ell$ and applying Lemma 8.8 we get

$$h(\Theta^m \tilde{\mu}, \sigma_p) \geq \log p - \frac{\log p}{\log p_*} \varepsilon.$$

Letting $\varepsilon \downarrow 0$ proves the theorem. □

## 9. Dimension of sum sets

For any measure $\mu$, let

$$(27) \qquad \dim \mu \stackrel{\text{def}}{=} \inf\{\dim_H S \mid S \text{ is a Borel set with } \mu(S) = 1\}.$$

Define the *lower dimension* of $\mu$ by

$$(28) \qquad \underline{\dim}\mu \stackrel{\text{def}}{=} \inf\{\dim_H S \mid S \text{ is a Borel set with } \mu(S) > 0\}.$$

To compute dimension of measures we use the following lemma from Billingsley [1].

BILLINGSLEY'S LEMMA. *Let $\mu$ be a positive finite measure on $\mathbb{T}$. Assume $K \subset \mathbb{T}$ is a Borel set satisfying $\mu(K) > 0$ and*

$$K \subset \left\{ x \in \mathbb{T}: \liminf_{\varepsilon \downarrow 0} \frac{\log \mu[B_\varepsilon(x)]}{\log \varepsilon} \leq \gamma \right\}.$$

*Then $\dim_H K \leq \gamma$. If the $\liminf$ is $\gamma$ a.e., then $\dim_H K = \gamma$.*

Here $B_\varepsilon(x)$ can be the interval of length $\varepsilon$ centered at $x$, or the mesh interval with edge $\varepsilon$ containing $x$, etc.

In most of this section we restrict attention to $p$-invariant measures on $\mathbb{T}$, for some fixed integer $p > 1$. In this case, by the Shannon-McMillan-Breiman (SMB) Theorem, it follows that if $\mu$ is ergodic then $\dim \mu = h(\mu, \sigma_p)/\log p$. If $\mu$ is not ergodic, denote by $\mu = \int \mu_\theta \, d\theta$ its ergodic decomposition, and then

(29) $$\dim \mu = \operatorname{ess\,sup}_\theta \dim \mu_\theta.$$

This (known) fact is proved in Meiri and Peres [15] in a more general context. We wish to derive the equivalent statement for lower dimension.

THEOREM 9.1. *Let $\mu$ be a $p$-invariant measure on $\mathbb{T}$, and denote its ergodic decomposition by $\mu = \int \mu_\theta \, d\theta$. Then $\underline{\dim}\mu = \operatorname{ess\,inf}_\theta \dim \mu_\theta$.*

*Proof.* Denote $\gamma \stackrel{\text{def}}{=} \operatorname{ess\,inf} \dim \mu_\theta$. By the above remarks,

$$\gamma = \operatorname{ess\,inf} h(\mu_\theta, \sigma_p)/\log p.$$

Let $\psi(x)$ denote the SMB function of $\mu$, i.e., $\psi(x) \stackrel{\text{def}}{=} \lim \frac{1}{n} I_{\alpha_0^{n-1}}(x)$ for the partition $\alpha = \{[\frac{j}{p}, \frac{j+1}{p}]\}_{j=0}^{p-1}$ (see Meiri and Peres [15] for more details). The function $\psi$ is constant on fibers, and the SMB theorem for nonergodic transformations implies that for almost every $\theta$ we have $\psi(x) = h(\mu_\theta, \sigma_p)$ $\mu_\theta$-a.e. (cf. Parry [16, p. 39]). Fix some $\varepsilon > 0$, and let

$$B = \left\{ x : \frac{\psi(x)}{\log p} < \gamma + \varepsilon \right\}.$$

From the definition of $\gamma$ we have $\mu(B) > 0$, hence $\underline{\dim}\mu \leq \dim_H B$. From Billingsley's lemma and the SMB theorem we know that $\dim_H B \leq \gamma + \varepsilon$. Since $\varepsilon$ was arbitrary, we conclude that $\underline{\dim}\mu \leq \gamma$. The other direction is proved similarly: if $\underline{\dim}\mu < \gamma$, then there exists a Borel set $B$ and some $\varepsilon > 0$ such that $\mu(B) > 0$ and $\dim_H B < \gamma - \varepsilon$. Applying once again Billingsley's lemma and the SMB theorem we get $\psi(x)/\log p \leq \gamma - \varepsilon$ for $x \in B$, contradicting the definition of $\gamma$. □

In particular, if $\mu$ is $p$-invariant and ergodic, we conclude that

$$\dim \mu = \underline{\dim}\mu = \frac{h(\mu, \sigma_p)}{\log p}. \tag{30}$$

LEMMA 9.2. *Any two measures $\mu$ and $\nu$ on $\mathbb{T}$ satisfy $\underline{\dim}(\mu * \nu) \geq \underline{\dim}\mu$.*

*Proof.* Suppose that $B$ is a Borel set with $\dim_H B < \underline{\dim}\mu$. We need to prove that $(\mu * \nu)(B) = 0$. For any $t \in \mathbb{T}$ we have $\dim_H(B - t) = \dim_H B < \underline{\dim}\mu$, so $\mu(B - t) = 0$. Hence $(\mu * \nu)(B) = \int \mu(B - t) \, d\nu(t) = 0$. □

COROLLARY 9.3. *Let $\mu$ and $\nu$ be $p$-invariant measures on $\mathbb{T}$, with $\mu$ ergodic. Let $\mu * \nu = \int \varphi_\theta \, d\theta$ be the the ergodic decomposition of $\mu * \nu$. Then $h(\varphi_\theta, \sigma_p) \geq h(\mu, \sigma_p)$ for almost-every $\theta$.*

*Proof.* Write

$$\operatorname*{ess\,inf} \frac{h(\varphi_\theta, \sigma_p)}{\log p} = \underline{\dim}(\mu * \nu) \quad \text{by Theorem 9.1}$$
$$\geq \underline{\dim}\mu \quad \text{by Lemma 9.2}$$
$$= \frac{h(\mu, \sigma_p)}{\log p} \quad \text{by (30)}.$$

This proves the assertion. □

The point in the last corollary is that $\mu * \nu$ need not, in general, be ergodic, even if $\mu$ and $\nu$ are ergodic. Furthermore, the entropy of its ergodic components need not be equal, as seen in the following example:

*Example* 9.4. A nonergodic convolution.

Take $p = 2$ and let $X = \sum_{i=1}^\infty x_i p^{-i}$ denote the random variable on $\mathbb{T}$ for which $x_i = 0$ if $i \not\equiv 0 \pmod 5$, and otherwise $x_i$ is 0 or 1 with probability $\frac{1}{2}$. Denote by $\mu$ the distribution of $X + pX + \cdots + p^4 X \pmod 1$. Then $(\mathbb{T}, \mu, \sigma_p)$ is ergodic but not weakly-mixing (since $\sigma_p^5$ is not ergodic). Also, $\mu^{*5}$ is not ergodic, and decomposes to (finitely many) components with different entropies (one of them is Lebesgue). It is also possible to construct an example with a continuum of components.

We turn now to topological corollaries of Theorems 1.1 and 1.8.

*Proof of Corollary* 1.2. By the variational principle for expansive maps (see Walters [19]), for every $i$ there exists a $p$-invariant and ergodic measure $\mu_i$ supported on $S_i$ with $h(\mu_i, \sigma_p) = h_{\text{top}}(S_i, \sigma_p)$. Recall that

$$h_{\text{top}}(S_i, \sigma_p)/\log p = \dim_H S_i;$$

hence $\dim \mu_i = h(\mu_i, \sigma_p)/\log p = \dim_H S_i$ by (30). By our assumptions, $h_i \stackrel{\text{def}}{=} h(\mu_i, \sigma_p)/\log p$ satisfy (1).

Denote $\nu^{(n)} = \mu_1 * \cdots * \mu_n$. Since $\nu^{(n)}(S_1 + \cdots + S_n) = 1$, we have
$$\dim \nu^{(n)} \leq \dim_H(S_1 + \cdots + S_n),$$
so it is enough to show that $\dim \nu^{(n)} \longrightarrow 1$. Denote by $\nu^{(n)} = \int \nu_\theta^{(n)} \, d\theta$ the ergodic decomposition of $\nu^{(n)}$. Then

$$\begin{aligned}
\dim \nu^{(n)} &= \operatorname{ess\,sup}\{\dim \nu_\theta^{(n)}\} && \text{by (29)} \\
&= \frac{1}{\log p} \operatorname{ess\,sup} h(\nu_\theta^{(n)}, \sigma_p) && \text{by (30)} \\
&\geq \frac{1}{\log p} \int h(\nu_\theta^{(n)}, \sigma_p) \, d\theta \\
&= \frac{1}{\log p} h(\nu^{(n)}, \sigma_p) && \text{by Rokhlin's theorem} \\
&\longrightarrow 1. && \text{by Theorem 1.1} \quad \square
\end{aligned}$$

*Proof of Theorem* 1.3. We define inductively a joining of full entropy $\tilde{\mu}$. Let $\nu^{(1)} = \mu_1$. Consider the ergodic decomposition of $\mu_1 \times \mu_2$. Since $\mu_1 \times \mu_2(S_1 \times S_2) > 0$, by Lemma 2.1 we can find an ergodic component $\nu^{(2)}$ such that

(i) $h(\nu^{(2)}, \sigma_p \times \sigma_p) = h(\mu_1, \sigma_p) + h(\mu_2, \sigma_p)$,

(ii) $\nu^{(2)}$ projects to $\mu_1$ and $\mu_2$, and

(iii) $\nu^{(2)}(S_1 \times S_2) > 0$.

Consider next the ergodic decomposition of $\nu^{(2)} \times \mu_3$, and find an ergodic component $\nu^{(3)}$ with similar properties. Continue in this manner to define a sequence of measures $\nu^{(n)}$ such that $\nu^{(n)}$ is an ergodic joining of full entropy of $\mu_1, \ldots, \mu_n$, and $\nu^{(n)}(S_1 \times \cdots \times S_n) > 0$. Let $\tilde{\mu}$ be the inverse limit of this system. Then, by Theorem 1.8, $h(\Theta^n \nu^{(n)}, \sigma_p) \longrightarrow \log p$. As $\nu^{(n)}$ is ergodic, $\Theta^n \nu^{(n)}$ is ergodic as well. Since $(\Theta^n \nu^{(n)})(S_1 + \cdots + S_n) > 0$, by (27) and (30) we conclude that

$$\dim_H(S_1 + \cdots + S_n) \geq \dim \Theta^n \nu^{(n)} = \frac{h(\Theta^n \nu^{(n)}, \sigma_p)}{\log p} \longrightarrow 1,$$

as stated. $\square$

## 10. Examples and questions

1. Recall that a measure $\mu$ on $\mathbb{T}$ is $\{c_n\}$-*normal* a.e. if $\{c_n x \pmod 1\}$ is uniformly-distributed for $\mu$-almost every $x \in \mathbb{T}$. Suppose that $\mu$ is a $p$-invariant ergodic measure, such that $\mu * \mu$ is $\{q^n\}$-normal. Does it follow that $\mu$ is $\{q^n\}$-normal as well? For noninvariant measures, the assertion is false:

*Example* 10.1.   $\mu * \mu$ normal does not imply that $\mu$ is normal.

Let $\{X_i\}$ and $\{Y_i\}$ be independent sequences of random bits with $\mathbb{P}(X_i = 0) = \frac{1}{2}$ and $\mathbb{P}(Y_i = 0) = 1/i$ for all $i \geq 1$. For $2^k \leq i < 2^{k+1}$ define $Z_i = X_i Y_k$. Let $\mu$ be the distribution of $\sum 2^{-i} Z_i$. Then $\mu$ is $\{2^n\}$-normal in probability, but not $\{2^n\}$-normal a.e. Also, $\mu * \mu$ *is* $\{2^n\}$-normal a.e.

2. Suppose that $\mu$ is a $p$-invariant and ergodic measure on $\mathbb{T}$ with positive entropy. Does it follow that $\underline{\dim}(\mu^{*n}) \to 1$ as $n \to \infty$? This would mean that the convergence $h(\mu^{*n}, \sigma_p) \to \log p$ is uniform on all ergodic components of $\mu^{*n}$.

3. As we remarked in the introduction, the entropy condition in Theorem 1.1 is sharp.

*Example* 10.2.   Given numbers $0 < h_i < 1$ with $\sum \dfrac{h_i}{|\log h_i|} < \infty$, we construct a sequence of $p$-invariant ergodic measures $\{\mu_i\}$ such that $h(\mu_i, \sigma_p)/\log p = h_i$, yet $\mu_1 * \cdots * \mu_n \not\to \lambda$ weak*.

For a sequence $\{\beta_i\}$, let $X_j^{(i)}$ be independent random variables on $\{0, \ldots, p-1\}$ with $\mathbb{P}(X_j^{(i)} = 0) = 1 - \beta_i$, and $\mathbb{P}(X_j^{(i)} = k) = \beta_i/(p-1)$ for $k = 1, \ldots, p-1$. Define $\mu_i$ to be the distribution of $\sum_{j=1}^{\infty} X_j^{(i)} p^{-j}$. Then $\{\mu_i\}$ are Bernoulli measures, and

$$h(\mu_i, \sigma_p) = H(1 - \beta_i, \frac{\beta_i}{p-1}, \ldots, \frac{\beta_i}{p-1}) \sim \beta_i \log \frac{1}{\beta_i}.$$

Define $\beta_i$ by requiring that $h(\mu_i, \sigma_p) = h_i \log p$. Then from the condition on $\{h_i\}$ we get $\sum \beta_i < \infty$ (see §6). It is not hard to see that $|\hat{\mu}_i(1) - 1| \leq 4\pi\beta_i$. Thus $\sum |\hat{\mu}_i(1) - 1| < \infty$. But then

$$|(\mu_1 * \cdots * \mu_n)^\wedge(1)| = |\prod_{i=1}^{n} \hat{\mu}_i(1)| \geq \prod_{i=1}^{\infty} |\hat{\mu}_i(1)| > 0,$$

whence $(\mu_1 * \cdots * \mu_n)^\wedge(1) \not\to 0$.

4. Is the dimension condition of Corollary 1.2 sharp as well? Specifically, given a sequence of numbers $0 < d_i < 1$ such that $\sum d_i/|\log d_i| < \infty$, can one always find $p$-invariant closed subsets $S_i \subset \mathbb{T}$ with $\dim_H S_i = d_i$ and $\lim \dim_H(S_1 + \cdots + S_n) < 1$? Currently we can construct sets satisfying the desired conclusion, but only when $\sum d_i/|\log d_i|$ is small enough.

*Example* 10.3.   Given numbers $0 < d_i < 1$ with $\sum \dfrac{d_i}{|\log d_i|} < c(p)$ (for some $c(p)$ that can be made explicit) there is a sequence of $p$-invariant closed sets $S_i$ such that $\dim_H S_i = d_i$, yet $\lim \dim_H(S_1 + \cdots + S_n) < 1$.

Indeed, consider the set

$$S(N, \beta) = \Big\{ x \in \mathbb{T} : \forall n \in \mathbb{N}, \quad (\# \text{ non-}0 \text{ digits in } x_{n\ldots n+N-1}) \leq \beta N \Big\}.$$

If $\mu(\beta)$ is constructed in the same way the measures $\mu_i$ were defined in the previous subsection, then clearly

$$\dim_H S(N, \beta) = \frac{h_{\text{top}} S(N, \beta)}{\log p} \longrightarrow \frac{h(\mu(\beta))}{\log p} \quad \text{as } N \longrightarrow \infty.$$

Thus we can chose $N_i$ and $\beta_i$ so that if $S_i = S(N_i, \beta_i)$, then

(31)
$$d_i < \dim_H S_i < (1+\varepsilon)d_i, \quad \text{and}$$
$$(1-\varepsilon)d_i < \frac{h(\mu(\beta_i))}{\log p} < (1+\varepsilon)d_i.$$

We also assume $N_i | N_{i+1}$ and $N_i \geq 2^i/\varepsilon$ for all $i$. The exact condition we need is that $\sum_{i=0}^{\infty} \beta_i < 1-\varepsilon$, which can be translated to a condition on $\sum_{i=0}^{\infty} d_i/|\log d_i|$ via (31).

We shall presently show that $\lim \dim_H(S_1 + \cdots + S_n) < 1$. It is possible to replace the sets $S_i$ by $S'_i \subset S_i$ such that $\dim_H S'_i = d_i$, using the fact that the $S_i$'s we defined above are shifts of finite type. Shifts of finite type have many closed $p$-invariant subsets, and in particular have a closed $p$-invariant subset with any Hausdorff dimension between 0 and the dimension of the full set. Thus there is a closed $p$-invariant subset $S'_i \subset S_i$ with the required Hausdorff dimension (properties of shifts of finite type are discussed in Denker et al. [5]; the above result is a consequence of the Jewett-Krieger theorem, (chap. 29 in that reference)). Clearly $\lim \dim_H(S'_1 + \cdots + S'_n) \leq \lim \dim_H(S_1 + \cdots + S_n)$. The proof that $\lim \dim_H(S_1 + \cdots + S_n) < 1$ follows from the following proposition:

PROPOSITION 10.4. *If $N|M$,*

$$S(N, \beta) + S(M, \beta') \subset S(M, \beta + \beta' + \tfrac{1}{M}).$$

*Proof.* Consider the $M$-block $z_{n\ldots n+M-1}$ for any $z \in S(N, \beta) + S(M, \beta')$ as an element of $\mathbb{Z}/p^M\mathbb{Z}$. Then there are $x \in S(N, \beta)$ and $x' \in S(M, \beta')$ such that

$$x_{n\ldots n+M-1} + x'_{n\ldots n+M-1} \bmod p^M \leq z_{n\ldots n+M-1}$$
$$\leq (x_{n\ldots n+M-1} + x'_{n\ldots n+M-1} \bmod p^M) + p - 1.$$

Note that for any $a, b \in \mathbb{N}$, their base-$p$ expansions satisfy
(32)
$$(\# \text{ non-}0 \text{ digits in } a) + (\# \text{ non-}0 \text{ digits in } b) \geq (\# \text{ non-}0 \text{ digits in } a+b).$$

This can be shown, for example, by induction.

Now in $x_{n\ldots n+M-1}$ there are at most $\beta M$ non-0 digits, and in $x'_{n\ldots n+M-1}$ at most $\beta' M$ non-0 digits. Thus using (32), $z_{n\ldots n+M-1}$ can have at most $(\beta + \beta')M + 1$ non-0 digits, and the proposition follows. □

Using the above proposition, we see that for any $n$,
$$S_1 + \ldots + S_n \subset S(N_n, \epsilon + \sum_{i=1}^{\infty} \beta_i).$$
As $\epsilon + \sum_{i=1}^{\infty} \beta_i < 1$, we see that
$$\lim_{n \to \infty} \dim_H(S_1 + \ldots + S_n) \leq \frac{h(\mu(\epsilon + \sum_{i=1}^{\infty} \beta_i))}{\log p} < 1.$$

5. Does Theorem A of Section 1.1 hold under the weaker assumption that the collision exponent is smaller than two? By the Bootstrap Lemma, a positive answer to the first question in this section would give a positive answer to this question. In particular, this would imply that every $p$-invariant measure of positive entropy is $\{q^n\}$-normal, for any $p, q$ such that some prime factor of $p$ does not divide $q$. The case where $\mu$ is a Bernoulli measure is covered by Feldman and Smorodinsky [7], for any multiplicatively-independent $p, q$.

6. The following example shows that Theorem 1.8 is not valid under the weaker assumptions of Theorem 1.1. In fact, even weak convergence is not guaranteed:

*Example* 10.5. A joining of full entropy $\tilde{\mu}$ with $\sum \frac{h_n}{|\log h_n|} = \infty$, yet $\Theta^n \tilde{\mu} \not\to \lambda$ weak$^*$.

Take $p = 2$, and fix some irrational $\alpha$. Let $t_0$ be a uniform random variable on the unit interval, and define $t_n = t_0 + n\alpha \pmod{1}$. Let $\{Y_i^j\}_{i,j=1}^{\infty}$ be a sequence of i.i.d. variables with $\mathbb{P}(Y_i^j = 0) = \mathbb{P}(Y_i^j = 1) = \frac{1}{2}$. Given a nonnegative sequence $\{h_j\}$, define
$$X_i^j = \begin{cases} 0 & \text{if } t_i \geq h_j, \\ Y_i^j & \text{if } t_i < h_j \end{cases}$$
and let $X^j = \sum_{i=1}^{\infty} X_i^j 2^{-i}$. Denote by $\mu_j$ the distribution of $X^j$, and by $\tilde{\mu}$ the distribution of $\{X^j\}_{j=1}^{\infty}$ on $\mathbb{T}^{\mathbb{N}}$. Then $\mu_j$ is a $\sigma_2$-invariant and ergodic measure with entropy $h_j$ (in base 2), and $\tilde{\mu}$ is a joining of full entropy. Take $h_j = \frac{1}{j+1}$. We claim that for small enough $\alpha$ we have $\Theta^k \tilde{\mu} \not\to \lambda$ weak$^*$. To see that, define the following events:
$$\begin{aligned} A_N &= \{\forall j \in \mathbb{N} \, \forall i = 1, \ldots, N \ X_i^j = 0\}, \\ B_N &= \{\forall i > N \, \forall j > 2^{i/2} \ X_i^j = 0\}. \end{aligned}$$
If $t_i \geq h_1$ then $X_i^j = 0$ for all $j$. Take $N$ such that $\sum_{i=N+1}^{\infty} 2^{-i/2} \leq \frac{1}{8}$, and fix $\alpha \in (0, \frac{1}{4N})$. Clearly $\mathbb{P}(A_N) \geq \mathbb{P}(t_i \geq h_1 \, \forall i = 1, \ldots, N) \geq \mathbb{P}(t_1 \in (\frac{1}{2}, \frac{3}{4})) = \frac{1}{4}$.

Also, $\mathbb{P}(B_N^c) \leq \mathbb{P}(\exists i > N, t_i < 2^{-i/2})$, since if $t_i \geq 2^{-i/2}$ for all $i > N$, then for all $j > 2^{i/2}$ we have $t_i > \frac{1}{j} > h_j$; hence $X_i^j = 0$. By our choice of $N$,

$$\mathbb{P}(B_N^c) \leq \sum_{i=N+1}^{\infty} \mathbb{P}(t_i < 2^{-i/2}) = \sum_{i=N+1}^{\infty} 2^{-i/2} \leq \frac{1}{8},$$

and we conclude that $\mathbb{P}(A_N \cap B_N) \geq \mathbb{P}(A_N) - \mathbb{P}(B_N^c) \geq \frac{1}{8}$. For $(X_i^j) \in A_N \cap B_N$, for any $k$ the first $N/4$ digits of the sum $X^1 + \cdots + X^k \pmod 1$ are zero, and so $\Theta^k \tilde{\mu} \not\longrightarrow \lambda$ weak$^*$.

*Acknowledgments.* Parts of this work are contained in the Ph.D. thesis of the second named author. We are grateful to Hillel Furstenberg, Dan Rudolph and Benjamin Weiss for helpful discussions. We also thank Bryna Kra, Russell Lyons and Wilhelm Schlag for useful comments on the presentation.

Institute of Mathematics, The Hebrew University, Jerusalem, Israel
*E-mail addresses*: elon@math.huji.ac.il
davidm@emc.com
peres@math.huji.ac.il